\documentstyle[10pt]{article}
\input amssym.def
\input amssym

\newtheorem{theorem}{Theorem}[section]
\newtheorem{lemma}[theorem]{Lemma}
\newtheorem{proposition}[theorem]{Proposition}
\newtheorem{corollary}[theorem]{Corollary}
\newtheorem{conjecture}[theorem]{Conjecture}
\newtheorem{definition}[theorem]{Definition}

\begin{document}
\newcommand{\Z}{\Bbb Z}
\newcommand{\R}{\Bbb R}
\newcommand{\Q}{\Bbb Q}
\newcommand{\C}{\Bbb C}

\begin{titlepage}
\title{Zagier's conjecture on $L(E,2)$}
\author{A.B. Goncharov,  A.M. Levin} 

\date{}
\end{titlepage}
\maketitle
\tableofcontents

\section{Introduction}

{\bf Summery}. In this paper we introduce an elliptic analog of the Bloch-Suslin complex and prove that it (essentially) computes the weight two parts of the groups $K_2(E)$ and $K_1(E)$
for an elliptic curve $E$ over an arbitrary field $k$.  Combining this with the results of Bloch and Beilinson we proved Zagier's conjecture on $L(E,2)$ for modular elliptic curves over $\Q$.

{\bf 1. The elliptic dilogarithm}. The dilogarithm is the following
multivalued analytic function of on $\Bbb CP^1 \backslash \{0,1,\infty\})$:
$$
Li_2(z) = -\int_0^z\log(1-t)\frac{dt}{t}
$$
It has a  single-valued version, the Bloch-Wigner function:
$$
{\cal L}_2(z):= Im Li_2(z) + \arg(1-z)\cdot \log|z|
$$

The elliptic analog of the dilogarithm was defined and studied by Spencer
Bloch in his seminal paper [Bl1]. 

The story goes as follows. Let $E(\Bbb  C) = \Bbb
C^{\ast}/q^{\Bbb Z}$ be the complex points
of an elliptic curve $E$. Here $q:= exp(2\pi i\tau ), Im
\tau >0$.

The function ${\cal L}_2(z)$ has a singularity of type $|z|\log|z|$ near
$z=0$. It satisfies the relation ${\cal L}_2(z) = -
{\cal L}_2(z^{-1})$. So averaging ${\cal L}_2(z)$ over the action of the group
$\Bbb Z$ on $\Bbb C^{\ast}$ generated by $z \longmapsto qz$ we get the
convergent series: 
$$
{\cal L}_{2,q}(z) := \sum_{n\in \Bbb Z} {\cal L}_2(q^nz),\qquad 
{\cal L}_{2,q}(z^{-1}) =  -{\cal L}_{2,q}(z)
$$
This function can be extended by linearity to the set of all divisors on
$E(\Bbb C )$ setting  ${\cal L}_{2,q}(P) :=
\sum_i n_i{\cal L}_{2,q}(P_i)$ for a divisor $P= \sum n_i (P_i)$. 

{\bf 2.  The results  on $L(E,2)$}. Let $L(E,s) =
L(h^1(E),s)$ be the Hasse-Weil 
$L$-function of an elliptic curve $E$ over  $\Bbb Q$.  
 We will always suppose that an elliptic curve $E$ has at least one
point over $\Bbb Q$: zero for the addition law.

Let  $v$ be a  
valuation  of a number field $K$, and $h_v$   the  corresponding canonical 
local  height on $E(K)$. As usual $x \sim_{\Bbb Q^{\ast}} y$
means that $x = qy$ for a 
certain  $q \in \Bbb Q^{\ast}$. Let $J = J(E)$ be the Jacobian of
$E$.

\begin{theorem} \label{zcc}
Let $E$ be a modular elliptic curve over $\Bbb Q$. Then there exists a
$\Bbb Q$-rational divisor $P = \sum n_j (P_j)$ over $\bar \Bbb Q$ which satisfy the conditions a)-c) listed below and such that 
\begin{equation} \label {resultaa}
L(E,2) \sim_{\Bbb Q^{\ast}} \pi \cdot {\cal L}_{2,q}(P)
\end{equation}
The conditions on   divisor $P$:
\begin{equation} \label {condition1}
a)  \qquad\qquad \qquad  \sum n_j P_j
\otimes P_j \otimes P_j =
0 \quad \mbox{in} \quad S^3J({\bar \Bbb Q}) \qquad\qquad\qquad
\end{equation}

b) For any  valuation $v$ of the field $\Q(P)$ 
generated by the coordinates of the points $P_j$ 
\begin{equation} \label {condition2}
 \sum n_j
h_v(P_j)\cdot P_j =0 \quad \mbox{in} \quad J({\bar \Bbb Q})\otimes \Bbb R
\end{equation} 

 c)  For every prime $p$  where $E$ has a   split multiplicative reduction 
one has  an integrality condition on $P$, see (\ref{icond121}) below.
\end{theorem}


 {\it The integrality condition }.    Suppose $E$ has a split multiplicative reduction at $p$ with
$N$-gon as  a special fibre. 
Let $L$ be a finite extention of $\Bbb Q_p$ of degree $n=ef$ and $ {\cal O}_L$ the ring of integers in $L$.
Let
$E^0$ be the connected component of the N\'eron model of $E$ over
$ {\cal O}_L$. 
Let us  fix an isomorphism 
$E^0_{F_{p^f}} = \Bbb G_m/{F_{p^f}}$. It
provides a
bijection between $\Bbb Z/{eN}\Bbb Z$ and the components of 
$E_{F_{p^f} }$. For a divisor $P$ such that all its points 
are defined over $L$
denote by $d(P;\nu)$   the degree of the restriction of the flat
extension of
a divisor $P$ to the $\nu$'th component of the $(eN)$-gon.

Let $B_3(x):= x^3 - \frac{3}{2}x^2 + \frac{1}{2}x$  be the third Bernoulli
polynomial.

 The integrality condition at $p$ is the following condition on a divisor $P$, provided by the work of Schappaher and Scholl ([SS]). For  a certain (and  hence
for any, see s. 3.3) extention
$L$ of $\Bbb Q_p$  such that
all   points of the divisor $P$ are defined  over $L$ one has ($[L:\Q_p]= ef$):
\begin{equation} \label {icond121}
\sum_{\nu \in \Bbb Z/(eN) \Bbb Z}d(P;\nu)B_3(\frac{\nu}{eN}) =0
\end{equation}

{\bf Remarks}. 1. For a $p$-adic valuation
$v$ of
the field $K(P)$ one has    $(\log p)^{-1} h_v(P_j) \in \Bbb Q$. So
the
condition b) in this case looks as follows
\begin{equation} \label {condition2nonar}
 (\log p)^{-1}\sum n_j
h_v(P_j)\cdot P_j =0 \quad \mbox{in} \quad J(K(P))\otimes \Bbb Q
\end{equation}
In particular the right hand side is a finite dimensional $\Bbb
Q$-vector space. 

2.  Lemma 1.5 below shows that, assuming  (\ref{condition1}),  if the condition (\ref{condition2}) is valid for all
archimedean valuations but one then it is valid for all of
them. In particular if $P \in \Bbb Z[E(\Bbb Q)]$   we can omit
(\ref{condition2}) for the archimedean valuation.

The proof of theorem
(\ref{zcc}) is based on 
the results of S. Bloch [Bl1] on regulators on elliptic 
curves,  a ``weak'' version of Beilinson's conjecture for  modular 
curves proved by A.A. Beilinson in [B2] and  the results presented in s.2-3
below.   

    To prove the theorem we introduce for an 
elliptic curve $E$ over an {\it arbitrary} field $k$ a new complex (the elliptic motivic complex $B(E;3)$) and prove  that its cohomology   essentially  computes   the weight 2 parts of $K_2(E)$ and $K_1(E)$   (see theorems (\ref{mrezz}) and  (\ref{zaza})). This complex  mirrors the properties of the elliptic dilogarithm.  
It is an elliptic deformation of
the famous Bloch-Suslin complex which computes $K_3^{ind}(F) \otimes \Bbb Q$  and  $K_2(F)$ for
an arbitrary field $F$ (see [DS], [S] and s.1.6). 
 
In particular we replace  the ``arithmetical'' condition b)  by its refined ``geometrical'' 
version (see s. 1.4), which   is equivalent to the condition b)
for curves over number fields. 
 
Our results imply   

\begin{theorem} \label {zccc}
Let $E$ be  an elliptic curve over $\Bbb Q$. Then

i) For any element $\gamma \in K_2(E)$ there exists a $\Bbb
Q$-rational divisor $P$ on $E$ satisfying the
conditions a), b) from theorem (\ref{zcc}) such that the value of
the Bloch-Beilinson regulator map $r_2: K_2(E) \longrightarrow
\Bbb R$ on $\gamma$ is $ \sim_{\Bbb
  Q^{\ast}}  {\cal L}_{2,q}(P)$

ii) For any $\Bbb
Q$-rational divisor $P$ on $E$ satisfying the
conditions a), b) there exists an element $\gamma \in
K_2(E)\otimes  \Bbb Q$ such that $r_2(\gamma)  \sim_{\Bbb
  Q^{\ast}} {\cal
  L}_{2,q}(P)$.
\end{theorem}

  Theorem (\ref{zccc}ii) implies immediately  

\begin{corollary} \label {zag11}
 Let $E$ be   an elliptic curve  over $\Q$. Let us assume   
that the image of $K_2(E)_{\Bbb Z} \otimes \Q$ under the regulator map is   $L(E,2) \cdot \Q$. (This is a part of the Bloch-Beilinson conjecture).

Then for any $\Bbb Q$-rational divisor $P $  on $E(\bar \Bbb
Q)$ satisfying the conditions a) -  c) of theorem (\ref{zcc})
one has
$$
q\cdot L(E,2) =  \pi \cdot {\cal L}_{2,q}(P)
$$
where $q$ is a rational number, perhaps equal to $0$.
\end{corollary}

{\bf Remark}.    Corollary (\ref{zag11}) has an analog 
for an elliptic curve over any number field.   Its formulation is an easy exercise to the reader. 

Unlike in Zagier's conjecture on $\zeta$-functions of
number fields one can not expect $P_i \in  E(\Bbb Q)$: the
Mordell-Weil group of an elliptic curve over $\Bbb Q$ could be trivial.

The conditions a)-b) are  obviously satisfied if $P$ is 
(a multiple of) a torsion divisor. Moreover, if $E$ is a curve with
complex multiplication then  $L(E,2)$ is  the value of the elliptic
dilogarithm on a torsion divisor ([Bl1]). However if $E$ is not a CM
curve this should not be true in general. Thus one has to  consider 
the  non-torsion divisors, and so it is
necessary to use the conditions a)-b) in full strength.

The 
conditions a) and b) were guessed by D.
Zagier several years ago after studying the results of the 
computer experiments with $\Bbb Q$-rational points on some elliptic
curves, which  he did with H. Cohen.

{\bf 3. A numerical example}. $E$ is given by equation $y^2 -y = x^3 -x$. The discriminant $\Delta$ $=$ conductor $= 37$. So $E$ has split multiplicative reduction at $p=37$ with one irreducible component of the fiber of the N\'eron model. Therefore the integrality condition is empty.

{\it Local nonarchimedean  heights on $E$}. Let $P=[ a/p^{2\delta}, b/p^{3\delta}] \in E(\Bbb  Q)$  where $a,b$ are prime to $p$.  If $p$ is prime  to $\Delta$  then $h_p(P) = 0$ if $\delta \leq 0$ and $h_p(P) = \delta \cdot \log p$ if $\delta > 0$. 
The local height at   $p=37$ is given by 
$h_{37}(P) = -1/6 + 2 \delta$ (see the formula  for the local height in  s. 4.3 of  or [Sil]).

The Mordell-Weil group has rank one and is generated by the point $P =[0,0]$. Consider the following integral points on $E$:
$$
P =[0,0], \quad 2P=[1,0], \quad 3P=[1,1],\quad  4P=[2,3],\quad  6P=[6,-14]   
$$
and also
$$
5P=[\frac{1}{4}, \frac{5}{8} ],\quad  10 P=[\frac{ 161}{16}, \frac{ 2065}{ 64}]
$$

 There  are no   height conditions  at $p \not = 37$ for the integral points and there is just one at $p=2$ for the points $5P$ and $10P$.

Consider the divisor $\sum n_k (kP)$. Notice that $S^3J(\Q) = \Z $ 
and  the   
 condition a) is $\sum n_k \cdot k^3 =0$.
The height condition at $p=37$  gives $\sum n_k \cdot k =0$ provided that the coordinates of $(kP)$ are   prime to $37$.

The divisor 
$$
P_k = (kP) - k(P) - \frac{k^3 -k}{6}((2P) - 2(P))
$$ 
satisfies the conditions $\sum n_k \cdot k = \sum n_k \cdot k^3 =0$. Also $P_{10} - 4 \cdot P_5$ satisfies the height condition at $p=2$.

The computer calculation (using PARI)  shows 
$$
 \frac{8 \pi \cdot{\cal L}_{2,q}( P_3)}{37 \cdot L(E,2) } = -8.0000...,  
 \qquad  \frac{8 \pi \cdot {\cal L}_{2,q}( P_4)}{37 \cdot L(E,2)  } = -26.0000..., 
$$
 $$
 \frac{8 \pi \cdot {\cal L}_{2,q}( P_6)}{37 \cdot L(E,2)  } = -90.0000..., \qquad
 \frac{8 \pi \cdot {\cal L}_{2,q}( P_{10} - 4 \cdot P_5)}{37 \cdot L(E,2)  } = -248.0000...
$$

{\bf 4. The group $B_2(E)$ and  a refined version of  conditions a) - b)}. Let $E$ be an elliptic curve
over an arbitrary field $k$ and  
$J:= J(k)$ be the group of $k$-points of the
Jacobian of $E$.  Let $\Bbb Z[X]$ be the free abelian group generated by a set $X$. 
We will define in s.  2.1   a group $B_2(E/k) = B_2(E)$ such  that  

a) one has 
an exact sequence 
\begin{equation} \label{exten}
0 \longrightarrow k^{\ast} 
 \longrightarrow  B_2(E/k) \stackrel{p}{\longrightarrow} S^2 J(k) \longrightarrow 0
\end{equation}

b) one has a canonical (up to a choice of a sixth root of unity) surjective homomorphism
\begin{equation} \label{homo}
h: \Bbb Z[E(k) \backslash 0] \longrightarrow B_2(E/k) 
\end{equation}
whose projection to $S^2J(k)$ is given by the formula 
$\{a\} \longmapsto a\cdot a$.

c) if $K$  is a  local field there is
a canonical homomorphism
$$
H: B_2(E/K) \longrightarrow \Bbb R
$$
whose restriction to the subgroup $K^{\ast} \subset B_2(E/K))$ 
is given by $x \longmapsto \log |x|$, ( see s. 2.3). 
Moreover the canonical local height $h_K$ is given by the composition
$$
\Bbb Z[E(K)\backslash 0] \stackrel{h}{\longrightarrow} B_2(E/K) \stackrel{H}{\longrightarrow} \Bbb R
$$

The group $B_2(E)$  appears naturally  as a version of  the theory of
biextensions.   It is  a ``motivic'' version of  
 theta functions. Set $\{a\}_2:= h(\{a\}) \in B_2(E)$.

The conditions a)-b) on a divisor $\sum_j n_j (P_j)$ are equivalent to the following single  one:
\begin{equation} \label {conditions1}
   \sum n_j
\{P_j\}_2 \otimes P_j =0 \quad \mbox{in} \quad B_2(E( {\bar \Bbb Q}))\otimes J({\bar \Bbb Q})
\end{equation}
 
More precisely,
\begin{lemma}
Let $K$ be a number field and $P_j \in E(K)$. Then
$$
\sum n_j
\{P_j\}_2 \otimes P_j =0 \quad \mbox{in} \quad B_2(E(K))\otimes
J(K)\otimes \Bbb Q
$$
if and only if the following two conditions hold:
$$
\sum n_j P_j
\otimes P_j \otimes P_j =
0 \quad \mbox{in} \quad S^3J(K) \otimes \Bbb Q 
$$
and for any valuation $v$ of the field $K$
$$
\sum n_j
h_v(P_j)\cdot P_j =0 \quad \mbox{in} \quad J(K)\otimes_{\Bbb Q} \Bbb R
$$
\end{lemma}

{\bf Proof}. Multiplying the exact sequence (\ref{exten}) by $J(K)\otimes \Bbb Q$ we get
$$
0 \longrightarrow K^* \otimes J(K)\otimes \Bbb Q \longrightarrow B_2(E(K))\otimes
J(K) \otimes \Bbb Q \stackrel{p \otimes id}{\longrightarrow} S^2J(K) \otimes J(K) \otimes \Bbb Q\longrightarrow 0
$$
and use the fact that the   local  norms $|\cdot |_v$ separate all the elements in $K^* \otimes \Bbb Q$.

{\bf 5. The elliptic motivic complex}. 
Let us suppose first that $k$ is an algebraically closed 
field. In  chapter 3 we define a subgroup $R_3(E) \subset \Bbb Z[E(k)]$. When $ k = \Bbb C$ it is a 
subgroup of all functional equations for the elliptic dilogarithm.
In particular the homomorphism
$$
{\cal L}_{2,q}: \Bbb Z[E(\Bbb C)] \longrightarrow \Bbb R, \quad \{a\} \longmapsto {\cal L}_{2,q}(a)
$$
annihilates the subgroup $R_3(E/\Bbb C)$.  
Consider the homomorphism ($J:=J(k)$)
$$
\delta_3: \Bbb Z[E(k)] \longrightarrow B_2(E) \otimes J, \quad \{a\} \longrightarrow
-\frac{1}{2} \{a\}_2 \otimes a
$$
  An important result (  theorem (\ref{pro})) is that  $\delta_3(R_3(E)) =0$ . Setting
$$
 B_3(E):= \frac{\Bbb Z[E(k)]}{R_3(E)}
$$
 we get  a homomorphism $\delta_3: B_3(E) \longrightarrow B_2(E) \otimes J$. 
Let us consider the following complex  
\begin{equation} \label {comp221}
B(E;3): \qquad B_3(E) \stackrel{\delta_3}{\longrightarrow}  B_2(E) \otimes J
\longrightarrow J \otimes \Lambda^2J \longrightarrow \Lambda^3J
\end{equation}
Here the middle arrow is $\{a\}_2 \otimes b \longmapsto a \otimes a\wedge
b$ and the last one is the canonical projection. The complex is
placed in degrees $[1,4]$. It is acyclic in the last two terms.  
 This is our 
elliptic motivic complex.

Let $I_E$ be the augmentaion ideal of the group algebra $\Z[E]$, and $I_E^4 $ its fourth pour.
Let $B_3^{\ast}(E)$  be the quotient of $I_E^4$ by the subgroup  generated by the elements $(f) \ast (1-f)^-$, where $\ast$ is the convolution in the group algebra $\Z[E]$, $f \in k(E)^*$,  and $g^-(t):= g(-t)$.   Then there is a homomorphism 
 \begin{equation} \label{pp}
\delta_3: B_3^*(E) \longrightarrow  k^* \otimes J
\end{equation}   
which  fits    the following     commutative diagram
$$
\begin{array}{ccccccc}
0&&0&&&&\\
\downarrow&&\downarrow&&&&\\
B_3^*(E)& \stackrel{ }{\longrightarrow} & k^* \otimes J &&&&\\
\downarrow &&\downarrow &&&&\\
B_3(E)&\stackrel{\delta_3}{\longrightarrow}  &B_2(E) \otimes J &\longrightarrow &J \otimes \Lambda^2J&\longrightarrow& \Lambda^3J\\
\downarrow &&\downarrow &&\downarrow = &&\downarrow =\\
S^3J &\longrightarrow &S^2J \otimes  J&\longrightarrow &J \otimes \Lambda^2J&\longrightarrow & \Lambda^3J\\
\downarrow&&\downarrow&&&&\\
0&&0&&&&
\end{array}
$$
 where the  vertical sequences are exact, and the bottom one is the Koszul  complex, and thus also exact.  
 Let us denote by $B^*(E;3)$ the complex (\ref{pp}). 
  It is canonically quasiisomorphic to the complex  $B(E;3)$.
 
The complex $B^*(E;3)$ looks simpler then   $B(E;3)$. However a definition of the differential in $B^*(E;3)$ which does not use the embedding to   $B(E;3)$ is   rather awkward, see s. 4.6.

If $k$ is not algebraically closed we 
postulate the Galois descent property:
$$
B(E/k;3) := B(E/{\bar k};3) ^{Gal(\bar k/k)}; \qquad  B^*(E/k,3) := B^*(E/{\bar k};3) ^{Gal(\bar k/k)}
$$

{\bf 6. Relation with algebraic $K$-theory}. Let $k$ be an arbitrary field. 
 Let ${\cal K}_2$ be  the sheaf of $K_2$ groups in the Zariski topology  on $E$. One has  canonical inclusion $K_2(k) \hookrightarrow H^0(E,{\cal K}_2)$ and   surjective  projection \linebreak $H^1(E,{\cal K}_2) \to  k^*$.
\begin{theorem}  \label{mrezz}   Let $k = \bar k$. Then there is a sequence
$$
  Tor(k^*, J) \hookrightarrow 
\frac{H^0(E,{\cal K}_2)}{K_2(k)} \longrightarrow  B_3^*(E) \longrightarrow  k^* \otimes J  \longrightarrow Ker (H^1(E,{\cal K}_2) \to k^*) \to 0
$$
    It is  exact in the term $k^* \otimes J$ and exact modulo $2$-torsion in the other terms.  
\end{theorem}

For an abelian group $A(E)$ depending functorially on $E$ let $A(E)^-$ be the  subgroup of  skewinvariants under the involution $x \to -x$ of $E$.   Recall that one has the  $\gamma$-filtration on the Quillen $K$-groups.
One can show that modulo $2$-torsion 
$$
\frac{H^0(E,{\cal K}_2)}{K_2(k)} = gr^{\gamma}_2K_2(E)^-, \qquad Ker (H^1(E,{\cal K}_2) \to k^*) = gr^{\gamma}_2K_1(E)^-
$$

Recall that the Bloch-Suslin complex for an arbitrary field $k$ is defined as follows:
$$
 B_2(k) \stackrel{\delta}{\longrightarrow} \Lambda^2k^*; \qquad B_2(k) := \frac{\Z [k^*]}{R_2(F)} ; \quad \delta: \{x\} \longmapsto (1-x) \wedge x
$$
Here $R_2(k)$ is the subgroup generated by the elements $\sum_i(-1)^i \{r(x_1,..., \hat x_i,... ,x_5)\}$, where $x_i$ runs through all $5$-tuples of distinct points over $k$ on the projective line  and $r$ is the cross ratio. 
 
 One should compare theorem (\ref{mrezz}) with the following exact sequence provided by  Suslin's  theorem on $K_3^{ind}(k)$  ([S]) and  Matsumoto's theorem on $K_2(k)$ ([M]) (see also a closely related results by Dupont and Sah [DS]):
$$
0 \longrightarrow \tilde Tor(k^*,k^*)  \longrightarrow K_3^{ind}(k) \longrightarrow B_2(k) \stackrel{\delta}{\longrightarrow} \Lambda^2k^* \longrightarrow K_2(k) \longrightarrow 0
$$ 
Here $ \tilde Tor(k^*,k^* ) $ is a nontrivial extension of   $Tor(k^*,k^* )$ by $\Z/2\Z$.
  
{\bf Remark}. To guess an elliptic analog of the Steinberg relation \linebreak $(1-x) \otimes x \in \C^* \otimes \C^*$ we might argue as follows. $E(\C) = \C^*/q^{\Z}$, so  let us try to make sence out of $\sum_{n \in \Z}(1-q^nx) \otimes q^nx$. Let $p:\C^* \to E(\C)$.   Projecting $\C^* \otimes \C^*$  to $\C^* \otimes J(\C) $ we 
get $\prod_{n \in \Z}(1-q^nx) \otimes x$, $x \in E(\C)$.
Regularizing the infinite product we   obtain $\theta(x) \otimes x$ where $$
\theta(x):= q^{1/12} z^{-\frac{1}{2}}
\prod_{j \geq  0}  (1-q^j z)\prod_{j  >  0} (1-q^j z^{-1})
$$
 Unfortunately this  seems to make no sence:  $\theta(x)$ is not a function on $E(\C)$ with values in $\C^*$,  but a section of a line bundle.  Only  introducing the {\it group} $B_2(E/\C)$ and realizing that $\theta(x)$ {\it is a function on $E(\C)$, but with values in  the group} $ B_2(E/\C)$, we   find  the elliptic analog of the Steinberg relation: $\theta(x)\otimes x \in  B_2(E/\C) \otimes J(\C)$ (compare with s. 4.1, 4.2, 4.6). For an arbitrary base field   $\{x\}_2 \in B_2(E)$ replaces $\theta(x)$.

In [W] J.Wildeshaus, assuming standard conjectures about mixed motives, gave a 
conjectural inductive definition of  groups similar to  $B_{n}(E)$ and
discuss an elliptic analog of weak version of Zagier's
conjecture.

{\bf Acknowledgement}.  This  paper was essentially  written  when the first author  enjoyed the hospitality of
 MPI(Bonn) and IHES during the Summer of 1995, and completed when the authors met in 1996 at the MPI (Bonn). The generous support of both institutions   gratefully
acknowledged. 

A.G. was 
partially supported by the  NSF Grant  DMS-9500010 and A.L. by
the 
Soros  International 
Scientific Foundation.

It is our pleasure to thank Don Zagier and Maxim Kontsevich for interesting 
discussions. We are very grateful to the referee  who   made a lot of useful remarks and pointed out some errors.

\section {The group $B_2(E)$}

{\bf 1. A construction of the group $B_2(X)$}. Let  
$k$  be an arbitrary field. 

For any two degree zero line bundles $L_1$ and $L_2$ on a regular curve 
$X$ over $k$ let us define, following Deligne ([De]), a $k^{\ast}$-torsor
$[L_1,L_2]$. 

{\it  Motivation}. Let $s_i$ be a section of a line bundle $L_i$. 
If  $div (s_2)=\sum m_i P_i$ then $<s_1,s_2>$ is
$s_1(div (s_2))\in {\otimes} L_1^{m_i}|_{p_i} = [L_1,L_2]$. 
Such a tensor product turns out to be 
symmetric and does not depend on the choice of $s_2$. If $L_1={\cal O}$ then $[L_1,L_2] =k^*$.

For $f \in k(X)^*$ and a closed point $x \in X_1$  set $\bar f(x):= Nm_{k(x)/k}f(x)$.  We extend $\bar f$ by linearity to the group of $k$-rational divisors on $E$.  

{\it Definition}. The  elements  of $[L_1,L_2]$ are pairs
$<s_1,s_2>$, where $s_i$ is a  section of the line bundle $L_i$ and the
divisors $div(s_1)$ and $div(s_2)$ are 
disjoint. 
 For a rational function $f,g$ such
that $div(f)$ is disjoint  from $div(s_2)$ and $div(g)$ is disjoint  from
$div(s_1)$,  one has ($(s):= div (s)$):
$$
<f\cdot s_1,s_2> =  \bar f((s_2)) <s_1,s_2>, \quad < s_1,g\cdot s_2> =  \bar  g((s_1)) <s_1,s_2>
$$
 There are two {\it a priory}   different  expansions
\begin{equation} \label{correct1}
<f\cdot s_1,g\cdot s_2> = \bar  f((g)) \bar f(s_2) < s_1,g\cdot s_2> = \bar  f((g))  \bar  f((s_2))   \bar g((s_1))<s_1,s_2>
\end{equation}
and
\begin{equation} \label{correct2}
<f\cdot s_1,g\cdot s_2> =  \bar  g((f)) \bar g((s_1)) < f\cdot s_1, s_2> =  \bar g((f)) \bar g((s_1)) \bar f((s_2))<s_1,s_2>
\end{equation}
The right hand sides coincide thanks to the Weil reciprocity $ \bar  f((g)) = \bar  g((f))$.
So the  $k^{\ast}$-torsor
$[L_1,L_2]$  is well defined. 

There is a canonical isomorphism of $k^{\ast}$-torsors
$$
[L_1 \otimes L_2, M ] \longrightarrow [L_1, M] \otimes_{k^{\ast}} [L_2, M ]
$$
and a similar additivity isomorphism for the second divisor. Further, 
the two possible natural isomorphisms
$$
[L_1 \otimes L_2, M_1 \otimes M_2 ] \longrightarrow  \otimes_{1 \leq i,j \leq 2 } [L_i , M_j] 
$$ 
coincide. 

Therefore  for any element $s \in S^2J_X$ one gets a $k^{\ast}$-torsor $[s]$
defined up to an isomorphism. Moreover $[s_1+s_2] = [s_1] \otimes_{k^{\ast}} 
[s_2]$. 

The  above facts just mean that the collection of $k^{\ast}$-torsors $[s]$
defines an   extension of type (\ref{exten}). This is the definition of the
group $B_2(X)$.

To construct the homomorphism (\ref{homo})   we want to make sense
of elements $<s_1,s_2>$ where the divisors of the sections $s_i$
{\it may not be
disjoint}. The definition is suggested by the motivation given above. 
Namely  
\begin{equation} \label{disj}
<s_1, s_2> \quad \in \quad [L_1,L_2] \otimes_{k^{\ast}} V(div(s_1),div(s_2))
\end{equation}
where $V(L_1,L_2)$ is a $k^{\ast}$-torsor defined as follows.
Recall that if $A \to B$ is a homomorphism of abelian groups  and $T$ is an $A$-torsor, we can define a $B$-torsor $f_*T := T \otimes _A B$ where $A$ acts on $B$ via the homomorphism $f$. 
Thus for a $k(x)^*$-torsor $T^*_xX$ we define a $k^*$-torsor $N(T^*_xX)$ using the norm homomorphism
$Nm: k(x)^* \to k^*$. 

For two arbitrary  divisors $l_1$ and $l_2$ on $X$ consider the 
following  $k^{\ast}$-torsor
$$
V(l_1,l_2):= \otimes_{x \in X_1}N(T^*_xX)^{\otimes ord_x l_1 \cdot ord_x l_2}
$$
Here $ord_x l$ is the multiplicity of the divisor $l$ at the point $x$. 
One has a canonical isomorphism 
$V(l_1 + l_2, m) \longrightarrow  V(l_1,m) \otimes V(l_2,m)$. 
If the divisors $l_1$ and $l_2$ are disjoint then $V(l_1,l_2) = k^{\ast}$.

Any rational function $f$ provides a canonical isomorphism
$$
 \varphi_f: V(l_1,l_2)  \longrightarrow  V((f) + l_1,l_2) \qquad v \longmapsto
{\tilde f}(l_2) \cdot v 
$$
In this formula   ${\tilde f}(x) \in N(T^*_xX)^{\otimes ord_x(f)}$ is  the ``leading
term'' of the function $f$ at   $x$. It is defined as follows.
Choose  a
local parameter  $t$ at the point $x$.  If $f(t) = at^k +$ higher order terms
 then ${\tilde  f}(x) = a(dt)^k \otimes 1 \in N(T^*_xX)^k = (T^*_xX)^k \otimes k^*$.  So
${\tilde  f(x)}^{ord_x l_2} \in N(T^*_xX)^{\otimes ord_x(f) \cdot ord_x l_2}$ and the
formula above makes sense.

Let us recall the full version of the
Weil reciprocity law.  Let
\begin{equation} \label{tame}
\partial_x:\{f,g\} \longmapsto (-1)^{ord_x(f)\cdot ord_x(g) } 
\frac{ f(x)^{ord_x(g)}}{ g(x)^{ord_x(f)}}
\end{equation}
be the tame symbol. Then 
for {\it any} two rational functions $f,g$ on a curve
over an algebraically closed field $k$ one has
$\prod_{x \in X} \partial_x(f,g) =1$.

There exists a canonical isomorphism of $k^*$-torsors
$$
S: [L_1,L_2] \otimes V((s_1),(s_2)) \longrightarrow [L_2,L_1] \otimes V((s_2),(s_1))
$$
 given on generators by the formula
$$
S: \quad <s_1,s_2> \quad \longmapsto\quad  (-1)^{deg L_1 \cdot deg L_2 + \sum_{x \in X_1} ord_xs_1 \cdot ord_xs_2}<s_2,s_1>
$$ 
 The defining properties of the torsor $[L_1,L_2]$ look as follows:
$$
<f\cdot s_1,s_2> = \bar f((s_2)) <s_1,s_2>, \quad < s_1,g\cdot s_2> = S(<s_1,(g)>) <s_1,s_2>
$$
The 
formula
\begin{equation} \label{oioioi}
S<(f),(g)> = <(f),(g)>  
\end{equation}
is equivalent to the Weil reciprocity. Similar to 
(\ref{correct1}), (\ref{correct2}) and using the formula (\ref{oioioi}) we see that $S$ is well defined.

{\bf 2.  The homomorphism (\ref{homo})}. From now on $X=E$ is an elliptic curve over $k$, so  there is canonical isomorphism $T_x^*E = T_0^*E$.  For   $c \in k^*$  one obviously has $<c\cdot s_1,s_2> = < s_1,s_2>$. Therefore we may  consider $< s_1,s_2>$ when $s_1,s_2$ are divisors on $E$. Consider a map
$$
\{a\}  \longmapsto <(a) - (0), (a) - (0)> \in  (T^{\ast}_0E)^{\otimes
2}\otimes_{k^{\ast}} [a\cdot a]
$$

There is an almost canonical (up to  a sixth root of unity) 
choice  of an element in $(T^{\ast}_0E )^{\otimes
2}$.  Namely, the quotient of $E$ by the involution $x \longmapsto -x$ is
isomorphic to $P^1$. The image of $0$ on $E$ is the 
point $\infty$ on $P^1$.  
The images of the three nonzero 2-torsion points on $E$
gives 3 distinguished points on $P^1$. 
Their ordering = choice of a level 2 structure on $E $.  
A choice of ordering  gives a {\it canonical} coordinate $t$ on 
$\Bbb A^1:=  P^1 \backslash \infty$ for which $t=0$ is the first point and $t=1$ is the second. 
 This coordinate provides a  vector in $T^{\ast}_{\infty}P^1$ and
so a vector in $(T^{\ast}_0E)^2$. Therefore we have six
different trivializations of $(T^{\ast}_0E)^2$ and thus a canonical one
of $(T^{\ast}_0E)^{12}$: their product. The sixth root of ($16 \times$  this
trivialization)  
is  the  (almost) canonical element in $(T^{\ast}_0E)^2$ we need.

Using this element in $(T^*_0E)^2$ we get a map (\ref{homo}). The composition
$\Bbb Z[E(k) \backslash 0] \longrightarrow 
B_2(E) \longrightarrow S^2J$  is obviously given by  $\{a\}  
\longmapsto a\cdot a$.

 If   $E$ is written in the
Weierstrass form $y^2 = (x-e_1) (x-e_2)(x-e_3)$  
then $e_1, e_2, e_3$ are the  coordinates of    the   distinguished points and 
$ \frac{x-e_i}{e_j -e_k}$ is the canonical coordinate corresponding to the ordering $e_i,e_j,e_k$.
  Let
$\Delta$ be the discriminant  of $E$. Then $\Delta =  16 \cdot \prod_{i<j}(e_i - e_j)^2$.
The canonical trivialization is $\Delta^{1/6}(dx/2y)^2$.

{\bf 3. The canonical height}. In this section we recall the
 construction of the canonical local heights via the
biextension (compare with  [Za] and  [Bl2]). The 
canonical local height gives a homomorphism $B_2(K)
\longrightarrow \Bbb Q$ (resp. to $\Bbb R$) when $K$ is a nonarchimedean
(resp. archimedean) local field.  

The construction of the group $B_2(E)$
provides us with a collection of $k^{*}$-torsors $T_{(x,y)}$ where
$(x,y)$ is a point of $J \times J$. The torsors $T_{(0,y)}$ and
$T_{(x,0)}$ are trivialized: we have a distinguished element 
$<\emptyset,(y) - (0)> \in T_{(0,y)}$ and a similar one in $T_{(x,0)}$.

The collection of torsors $T_{(x,y)}$ glue to a $k^{*}$-bundle on
$J \times J$. It is isomorphic to the (rigidified) Poincar\'e line bundle minus
zero section. 
 
Let $L$ be a degree zero line bundle
on an  elliptic curve over $k$. Let us denote by $L^{*}$ the
complement  of the zero section in $L$. It is a principal
homogeneous space over a certain commutative algebraic group
$A(L)$ over $k$ which is an extension
$$
0 \longrightarrow \Bbb G_m \longrightarrow  A(L) \longrightarrow
J \longrightarrow 0 
$$
The group $A(L)$ is described as follows. 
For any $a \in E$ let $t_a: x \longmapsto a+x$ be the
shift  by $a$. Then $t_a^{*} L$ is isomorphic to $L$ (because $L$ is of
degree zero). The set of isomorphisms from $L$ to $t_a^{*} L $
form a $k^{*}$-torsor. These torsors glue together to a
$k^{*}$-torsor over $E$ which is isomorphic to $L^{*}$. 
On the other hand the collection of isomorphisms  $
L \longrightarrow  t_a^{*} L$ form a group. This is the
group $A(L)$. It is
commutative: the commutator provides a morphism from $E
\times E$ to $\Bbb G_m$, which  has to be  a constant map. 

Now let $K$ be a local field and $E$ be an elliptic curve over $K$. We get 
a group extension
$$
0 \longrightarrow K^{\ast} \longrightarrow  A(L)(K) \longrightarrow J(K)\longrightarrow 0
$$
Let $U(G)$ be the maximal compact subgroup of a locally compact commutative
group $G$. One can show (use lemma 6.1, ch. 11 in [La]) that
$U\Bigl(A(L)(K)\Bigr)$ projects surjectively onto $J(K)$ if $K$ is 
archimedean or $E$ has a good reduction over
$K$. In the case of bad reduction the image is a subgroup of finite index.

There is canonical homomorphism 
$$
A(L)(K) \longrightarrow A(L)(K)/U\Bigl(A(L)(K)\Bigr) =:H
$$
The
quotient $H$ is isomorphic to $\Bbb Z$ (resp. to a subgroup in
$\Bbb Q$ which is an extension of 
$\Bbb Z$ by a finite group ) when $K$ is nonarchimedean and 
$E$ has good reduction (resp. bad reduction), and to 
$\Bbb R$ if $K= \Bbb C, \Bbb R$. 

Therefore we get a homomorphism
$
A(L)(K) \longrightarrow \Bbb Z (\mbox{resp}\quad \Bbb Q )
$
for the nonarchimedean case and 
$
A(L)(K) \longrightarrow \Bbb R 
$
for the archimedean one.

For a given $x$ the torsors $T_{(x,y)}$ form a group
$T_{(x,\cdot)}$ which is isomorphic to the group $A(L|_{x \times
  J})$, and there is a similar statement for the torsors  $T_{(x,y)}$
for a given $y$. Applying the homomorphism 
$
A(L|_{x \times J}) \longrightarrow H
$ 
we will get a homomorphism of the group $T_{(x,\cdot)}$ to 
$H$. Similarly we have a homomorphism of the group $T_{(\cdot, y)}$
to $H$. 

Consider the map
$$
U(T_{(x_1,\cdot)}) \times U(T_{(x_2,\cdot)}) \longrightarrow
T_{(x_1+x_2,\cdot)} 
$$
induced by the multiplication $T_{(x_1,y)}\times
T_{(x_2,y)} \longrightarrow
T_{(x_1+x_2,y)}$. Its image is a subgroup. This follows from the
commutativity of the diagram
$$
\begin{array}{ccc}
T_{(x_1,y_1)} \times T_{(x_2,y_1)}\times T_{(x_1,y_2)} \times
T_{(x_2,y_2)}&\longrightarrow& T_{(x_1,y_1+y_2)} \times T_{(x_2,y_1+y_2)}\\
&&\\
\downarrow&&\downarrow\\
&&\\
T_{(x_1+x_2,y_1)} \times T_{(x_1+x_2,y_2)}&\longrightarrow&T_{(x_1+x_2,y_1+y_2)}
\end{array}
$$
It is therefore a compact subgroup, and so a maximal compact
subgroup in $T_{(x_1+x_2,\cdot)}$. 
In particular  
the restrictions of the homomorphisms $T_{(x,\cdot)} \longrightarrow H$
and $T_{(\cdot,y)} \longrightarrow H$ to  $T_{(x, y)}$
coincide. 

So we get a well defined homomorphism 
$B_2(E(K)) \longrightarrow H$.  
Now the composition 
\begin{equation} \label{heig}
\Bbb Z[E(K)\backslash 0] \stackrel{h}{\longrightarrow}
B_2(E(K)) \longrightarrow H 
\end{equation}                               
is  
the canonical N\'eron height on $E(K)$. The restriction of the
homomorphism (\ref{heig}) to the subgroup $K^{*} \subset B_2(E(K))$
coincides with the logarithm of the norm homomorphism. In the 
nonarchimedean case $H$ is a subgroup of $\log p \cdot \Bbb Q
\subset \Bbb R$. 

{\bf Remark}. The homomorphism $h : \Bbb Z[E(K)\backslash 0]
\longrightarrow B_2(E(K))$ was defined up to a sixth root of
unity. This does not affect the definition of the height because
the norm vanishes on roots of unity.

\section{An elliptic analog of the Bloch-Suslin complex}

{\bf 1. The group $B_3(E)$ and complex $B(E;3)$}. We will assume in s.3.1 - 3.3 that $k = \bar k$. 
We will always use notation $J := J(k)$. Set $g^-(t) := g(-t)$. Denote by $\ast$ the convolution in the group algebra $\Bbb Z[E(k)]$. 
\begin{definition}
 $R_3(E)$  is the subgroup of $\Bbb Z[E(k)]$ 
generated by the elements 
$
(f) \ast (1-f)^-, \quad f \in k(E)^{\ast}
$, $\{0\}$, 
  and  the "distribution relations"
$$
m\cdot(\{a\} - m\cdot \sum_{mb=a}\{b\}), \quad a \in E(k), \quad m = -1,2
$$
\end{definition} 
 
 {\bf Remarks}. 
a) For $m=-1$ we  get the elements $\{a\} + \{-a\} \in R_3(E)$.  If we remove them from the definition of $R_3(E)$, we   get the same group. 

b) It would be more natural to add to the subgroup $R_3(E)$  
 the distribution relations
 for all $m \in \Bbb Z \backslash 0$: 
 we should  get the same group (compare with   lemma (\ref{destr})). But we will not need this.

Consider the  homomorphism
$$
\beta: \otimes^2k(E)^{\ast} \longrightarrow \Bbb Z[E(k)]
\qquad \beta: f\otimes g \longmapsto  f*g^- := \sum n_im_j\{a_i - b_j\}
$$
(the Bloch map), where $(f) = \sum n_i (a_i)$ and $(g) = \sum m_j(b_j)$.

 Recall that 
$$
\delta_3: \Bbb Z[E(k)] \longrightarrow   B_2(E) \otimes J, \qquad
\{a\} \longmapsto -\frac{1}{2} \{a\}_2 \otimes a
$$
and  $i: k^{\ast} \hookrightarrow B_2(E)$ is the canonical embedding
(see (\ref{exten})).  
Let   $I_E$ be the augmentation ideal of the group
algebra $\Bbb Z[E(k)]$ and   $p: I_E \to J$   the canonical projection.

Recall that if $k = \bar k$ the tame symbol provides a homomorphism
$$
\otimes^2 k(E)^* \stackrel{\partial}{\longrightarrow}   k^{\ast} \otimes \Z[E] 
$$
The Weil reciprocity law shows that its image belong to 
$k^{\ast} \otimes I_E$.
\begin{theorem} \label{commm}
The following diagram is commutative   
$$
\begin{array}{ccc}
\otimes^2 k(E)^* &\stackrel{\partial}{\longrightarrow}&  k^{\ast} \otimes I_E \\
&&\\
\downarrow \beta&&\downarrow i \otimes p\\
&&\\
\Bbb Z[E(k)]&\stackrel{\delta_3}{\longrightarrow}& B_2(E) \otimes J 
\end{array}
$$
\end{theorem}
 
{\bf Proof}. Let $
(f) = \sum n_i(a_i),  \quad (g) =  \sum m_j (b_j)$. 
Then
$$
\delta_3 \circ \beta (f\wedge g) = -\frac{1}{2} \cdot \sum_{i,j} m_in_j<a_i -b_j, a_i -b_j>
\otimes (a_i -b_j)
$$
The term 
$
\sum_{i,j} m_i n_j<(a_i -b_j), (a_i -b_j)>
\otimes a_i 
$
 equals to 
$$
\sum_{i,j} m_i n_j<(a_i) -(0), (a_i) -(0)>
\otimes a_i + \sum_{i,j} m_i n_j<(b_j) -(0), (b_j) -(0)>
\otimes a_i - 
$$
$$
2\cdot \sum_{i,j} m_i n_j<(a_i) -(0), (b_j) -(0)>
\otimes a_i
$$
The first term here is  zero because $\sum_{j} m_j =0$. The second is
zero because $\sum_{i} n_ia_i =0$ in $J$. The last one can be written as
$$
-2 \cdot \sum_{i} m_i <(a_i) -(0), (g)>
\otimes a_i 
$$
So the theorem follows from   the definition of the tame
symbol. 

\begin{theorem} \label{pro}
$\delta_3(R_3(E)) =0$. 
\end{theorem}

{\bf Proof}.   We will denote by $\{a\}_3$   projection of the generator $\{a\}$ onto the quotient $B_3(E)$. The map  $\delta_3$ kills the distribution relations:
$$
\delta_3 \Bigl(m(\{a\}_3 - m\sum_{mb=a}\{b\}_3))\Bigr) = m (\{a\}_2 \otimes a 
- m\sum_{mb=a}\{b\}_2 \otimes b) =
$$
$$
 m(\{a\}_2 - \sum_{mb=a}\{b\}_2)\otimes a =0
$$
The last equality is provided by  corollary (\ref{divr}), which will be proved later, in s. 4.4. (The proof   does not depend on any results or constructions in chapter 3).
 
The fact that   $\delta_3 (f \ast (1-f)) =0$ lies deeper and 
follows  from the theorem (\ref{commm}). Theorem (\ref{pro}) is proved.

Set
$B_3(E):= \Bbb Z[E(k)]/R_3(E)$.   
We get a complex 
$$
B(E,3):= \qquad B_3(E) \longrightarrow   B_2(E) \otimes J
\longrightarrow J \otimes \Lambda^2J \longrightarrow \Lambda^3J
$$
  

 Let  $I^k_E$  be $k$-th  power of the augmentation ideal. 
\begin{lemma} \label{cor}
$\beta( f\otimes g) \in I^4_E$. Moreover, $\beta$ is  surjective onto $I^4_E$.
\end{lemma}

{\bf Proof}. A divisor $\quad \sum n_i (a_i) \quad$ is principal if and only if $\quad \sum n_i = 0 $   and 
 $\sum n_i a_i =0$  in $J(E)$. So  $I_E^2$ coincides with the subgroup of $\Bbb Z[E]$
 given by the divisors of functions. So the convolution of two principal 
divisors belongs to $I^4_E$ and, moreover,  generate it. 

Let $B_3^{\ast}(E)$  be the  quotient of $I_E^4$  by the subgroup generated by the elements $(f) \ast (1-f)^-$. 
 
\begin{lemma} \label{tl}
$\delta_3(I^4_E) \in k^{\ast}\otimes J
 \subset B_2(E)\otimes J$
\end{lemma}
 
{\bf Proof}. It is easy to see that $\delta_3(I^4_E)  \subset h(I_E^3) \otimes J$. Further, $h(I_E^3) \subset k^*$ because $p \circ h (I^3_E)=0$ (see the  properties of the group $B_2(E)$ listed in s. 1.4). The lemma follows.

 So we get  a  complex 
\begin{equation} \label {comp2221}
B^*(E;3): \quad B_3^{\ast}(E) \stackrel{\delta_3}{\longrightarrow}  k^{\ast} \otimes J
\end{equation}

{\bf 2. Relation with algebraic $K$-theory}.  Let us  remind the long exact sequence of localization    
$$
K_3(k(E)) \stackrel{\partial_3}{\longrightarrow}   \oplus_{x \in E_1} K_2(k(x))    \longrightarrow  K_2(E) \longrightarrow      
$$
$$
K_2(k(E)) \stackrel{\partial_2}{\longrightarrow}  \oplus_{x \in E_1} k(x)^*   \longrightarrow  K_1(E) \longrightarrow  k(E)^*   
\stackrel{\partial_1}{\longrightarrow}  \oplus_{x \in E_1} \Bbb Z    
$$
 The group $K_1(E)$  has a subgroup $k^*$ which comes from the base.  One can show that
$$
H^0(E, {\cal K}_2) = gr^{\gamma}_2K_2(E), \qquad H^1(E, {\cal K}_2) = gr^{\gamma}_2K_1(E) = K_1(E)/k^*  
$$ 
   \begin{lemma} \label{3.7} Modulo $2$-torsion  one has
$$
  H^0(E, {\cal K}_2)  = K_2(k) \oplus H^0(E, {\cal K}_2)^-
$$
$$
 Ker \Bigl( H^1(E,{\cal K}_2) \longrightarrow k^* \Bigr)  =  H^1(E,{\cal K}_2)^-
 = K_1(E)^-
$$ 
\end{lemma}
 
 {\bf Proof}. It follows easily using the transfer related to the projection $E \to P^1$ given by factorization along the involution $x \to -x$. 

Recall that $H^3(B(E,3)) = H^4(B(E,3)) = 0$.
\begin{theorem} \label{za}
 Let $k = \bar k$. Then the commutative diagram from theorem (\ref{commm})  provides a  morphism of complexes
$$ 
\begin{array}{ccc}
\frac{K_2(k(E))}{k^{\ast} \cdot k(E)^{\ast}}
&\stackrel{\partial}{\longrightarrow}&  k^{\ast} \otimes J \\ 
&&\\
\tilde  \beta \downarrow&&\downarrow id\\ 
&&\\
B^*_3(E)&\stackrel{\delta_3}{\longrightarrow}&  k^* \otimes J
\end{array}
$$
where $ \tilde \beta$ is surjective and $Ker  \tilde \beta = Tor(k^*,J)$ modulo $2$-torsion.
\end{theorem}

This theorem and lemma (\ref{3.7}) implies immediately
\begin{theorem} \label{zaza}
   Let $k = \bar k$. Then  there   are an embedding 
$$
i: Tor(k^*, J(k)) 
\hookrightarrow \frac{H^0(E, {\cal K}_2)}{ K_2(k)}
$$
 and   canonical isomorphisms 
\begin{equation} \label{098}
  \frac{H^0(E, {\cal K}_2)}{ Tor(k^*, J(k))  +K_2(k)} =  H^1B^*(E;3)  \quad \mbox{modulo $2$-torsion}
\end{equation} 
\begin{equation} \label{099} 
Ker \Bigl( H^1(E,{\cal K}_2) \longrightarrow k^* \Bigr) = H^2B^*(E;3)
\end{equation} 
For an arbitrary field $k$ (\ref{098}) and (\ref{099}) are   isomorphisms modulo torsion.  
\end{theorem}

{\bf Proof of theorem ({\ref{za})}}.  We will first prove that we have a morphism of complexes and $Tor(k^*,J) \subset Ker \tilde \beta$, and then that $Tor(k^*,J) = Ker  \tilde \beta$ modulo $2$-torsion. 

    Recall that  the tame symbol homomorphism $\partial$ maps $K_2(k(E))$ to  $k^{\ast} \otimes  I_E$. Therefore the complex  computing the groups $H^0(E, {\cal K}_2)$ and $Ker \Bigl( H^1(E,{\cal K}_2) \longrightarrow k^* \Bigr)$  looks as follows
$$
K_2(k(E)) \stackrel{\partial}{\longrightarrow}  k^{\ast} \otimes I_E 
$$
Notice that $\partial(\{c,f\}) = c \otimes (f)$, so  it has a subcomplex 
 \begin{equation} \label{lllll}
k^{\ast} \cdot k(E)^{\ast}\stackrel{\partial}{\longrightarrow}  k^{\ast} \otimes  I_E^2
\end{equation} 
where $\cdot$ is the product in $K$-theory.  One obviously has $Ker \partial = K_2(k)$, and factorising by $K_2(k)$ we get the identity map $k^{\ast} \otimes  I_E^2 \to k^{\ast} \otimes  I_E^2$.

 The map $\beta: \otimes^2  k(E)^* \to I^4_E$ followed by the natural projection \linebreak $I^4_E \to  B^*_3(E)$ leads to  surjective map
$$
 \tilde \beta:   \otimes^2  k(E)^* \longrightarrow B^*_3(E)
$$
Notice that $\beta(k^{\ast} \otimes k(E)^{\ast}) =0$ and $\beta(f \otimes (1-f)) =0$ by the definition of the group $B^*_3(E)$. 
 So we get the desired morphism of complexes.

Tensoring the exact sequence $0 \to I^2_E \to I_E \to J \to 0$ by $k^*$ and using the fact that $I_E$ is a free abelian group, we get an exact sequence
\begin{equation}   \label{eexxs}
0 \longrightarrow Tor(k^*,J) \longrightarrow k^* \otimes I_E^2\longrightarrow k^* \otimes I_E \stackrel{id \otimes p}{\longrightarrow} k^* \otimes  J\longrightarrow 0
\end{equation}
So we get the following commutative diagram, where the vertical sequences are   complexes,   the complex on  the right  is the exact sequence (\ref{eexxs}), and $\alpha $ is injective:
$$
\begin{array}{ccc}
&&0\\
&&\downarrow \\
  &&Tor(k^*,J)\\
  &\swarrow&\downarrow \\
k^* \otimes I_E^2 &\stackrel{=}{\longrightarrow}&k^* \otimes I_E^2\\
 \alpha \downarrow&&\downarrow\\
 \frac{K_2(k(E))}{K_2(k)}&\stackrel{\partial }{\longrightarrow}&k^* \otimes I_E\\
&&\\
\tilde  \beta \downarrow&&\downarrow \\
 B_3^*(E)&\stackrel{\delta_3 }{\longrightarrow}&k^* \otimes J\\
  &&\downarrow\\
 &&0
\end{array}
$$ 
From this diagram we see that 
$
Tor(k^*,J) \subset Ker \tilde \beta   
$. 
Notice that
$$
Coker \Bigl(k^* \otimes I_E^2 \stackrel{\alpha}{\longrightarrow}  \frac{K_2(k(E))}{K_2(k)}\Bigr) = \frac{\otimes^2I^2_E}{\{(1-f) \ast(f)^-\}}
$$

  \begin{theorem} \label{mmmm}
 The map $\bar \beta: \frac{\otimes^2I^2_E}{\{(1-f) \ast(f)^-\}}\longrightarrow
B^{\ast}_3(E)$ is an isomorphism modulo $2$-torsion. 
\end{theorem}

{\bf Proof of theorem (\ref{mmmm})}. It consists of three steps of quite different nature. 

{\it Step 1}.  Set  $t_a: x \to x+a$, $(t_af)(x):=
f(x-a)$.

\begin{proposition}  \label{ref} 
Let $f$ and $g$ be rational functions
on   $E$. Then  
$$\{f,g\} - \{t_af,t_ag\} =0 \quad {\rm  in } \quad 
\frac{ K_2(k(E))}{k^*\cdot (k(E))^*}$$
\end{proposition}

Let $L/K$ be an extension of fields. Then one has a natural map $p^*: K_2(K) \to K_2(L)$ and the transfer map $p_* : K_2(L) \to K_2(K)$. 
  We need the following result ([BT]).

\begin{lemma}   \label{BT}
Let $L/K$ be a degree 2  extension of fields.  Then
 $K_2(L)$ is generated by symbols $\{k,l\}$ with $k \in K, l \in L$ and $p_*(\{k,l\}) = \{k,Nm_{L/K}l\}$. 
\end{lemma}

 In particular 
$p^*p_* =  Id + \sigma $  where $\sigma$ is a nontrivial element in  $Gal(L/K)$, 
 and thus  modulo $2$-torsion any Galois invariant element in $K_2(L)$  belongs to    $p^*K_2(K)$.

\begin{lemma}  \label {4.11}
Assume $k  = \bar k$.
Then any  rational 
function $f$ on an elliptic curve 
$E$ over $k$ can be decomposed into  a product
of   functions with divisors of the
following kind: $(a)-2(b)+(-a+2b)$.
 \end{lemma}

{\bf Proof of the lemma \ref{4.11}}. Any  
function $f$   can be decomposed into a   product
of the functions with divisors of the
following kind: $(a)-( b)-( c)+(-a+b+c)$.

Indeed,  let
$(f)=\sum n_a (a)$. We will use induction on $\sum|n_a|$. 
 Since  $\sum n_a =0$,   replacing if needed $f$ by $f^{-1}$  we can  find  
points $a,b,c $ such that
   $n_a>0$,   $n_b < 0$, $n_c < 0$ . Then  $ (f)-[(a)-(b)-(c)+(-a+b+c)]$ is    a  principal divisor  with smaller
 $\sum|n_a|$.  

 Further  $(a)-(b)-(c)+(-a+b+c)=
[(a)-2(d)+(-a+b+c)]- [(b)-2(d)+(c)$ for
$2d=b+c$.

{\bf Proof of the  proposition (\ref{ref})}.  According to the lemma we can assume that  $(f)=(b)-2(0)+(-b)$, $(g)=(c)-2(d)+(-c+2d)$.
Therefore 
$(t_af)=(b+a)-2(a)+(-b+a)$, $(t_ag)=(c+a)-2(d+a)+(-c+2d+a)$.

The quotient of $E$ under the involution $\sigma _a: x \to a-x$ is  isomorphic 
to $\Bbb P^1$. The symbol
$\{t_af,t_ag\}+\{\sigma _at_af,\sigma _at_ag\}$
is $\sigma _a$-invariant, so it comes from
   $K_2(k(\Bbb P^1) )$. It is known
that $K_2 (k(\Bbb P^1) )$ is generated
by $k(\Bbb P^1)^*\otimes k^*$.  
 So 
$$\{f,g\}-\{t_af,t_ag\} \sim \{f,g\}+\{\sigma _at_af,\sigma _at_ag\}=
\{f,g\cdot \sigma _0g\}\sim 0$$
where $x \sim y$ means  $x-y =0$ in $\frac{ K_2(k(E))}{k^*\cdot (k(E))^*}$. Indeed,  $\sigma _at_a =\sigma _0$, $\sigma _0f=f$,
and the symbol $\{f,g\cdot \sigma _0g\}$ is $ \sigma_0 $-invariant. Therefore it is $\sim 0$ by lemma (\ref{BT}).

{\it Step 2}. Let $A$ be an abelian group.   Let $S_k(A) \subset S^kI_A $ be the subgroup generated by the elements 
$$
 ((x_1)\ast X_1) \circ ... \circ  ((x_k) \ast X_k)  - X_1 \circ ... \circ  X_k, \quad x_1 + ... + x_k =0, \quad x_i \in A, \quad  X_i \in I_A
$$  
  This subgroup clearly belongs to the kernel of the convolution  map
$$
S^kI_A \to I_A^k, \quad X_1 \circ ... \circ X_k \longmapsto X_1 \ast ... \ast X_k
$$
So we get a homomorphism 
$
\alpha_k: S^kI_A/S_k(A) \longrightarrow I_A^k
$. 
\begin{proposition}  \label{schift}
 For any abelian group $A$  the homomorphism $\alpha_k$ is injective.
\end{proposition}

{\bf Proof}. We may assume that $A$ is finitely generated. Therefore 
 the   group ring $\Z[A]$ looks as follows:
\begin{equation} \label{Lorant} 
\Z[A] =   \Z[t_1,...,t_a, t_1^{-1},...,t_a^{-1}] \times \prod_{a+1 \leq j \leq   m}\frac{\Z[t_j]}{(t_j^{N_j} -1)}
\end{equation} 
 Under this isomorphism the augmentation ideal $I_A$ goes to the maximal ideal 
  $(t_1  - 1,...,t_{m}-1)$.
The subgroup $S_k(A)$ is generated by the elements
 \begin{equation} \label{rant}
(t_{i_1}  - 1) f_1 \circ ...  \circ (t_{i_{k-1}}  - 1) f_{k-1}   \circ (t_{i_k}  - 1) f_k \quad -  
\end{equation}
$$
(t_{i_1}  - 1)   \circ  ... \circ (t_{i_{k-1}}  - 1) \circ (t_{i_k}  - 1) f_1   ...   f_k
$$ 
So any element of the quotient $S^kI_A / S_k(A)$ can be written as
\begin{equation} \label{orant}
 (t_{i_1}  - 1)   \circ  ... \circ (t_{i_{k-1}}  - 1) \circ (t_{i_k}  - 1) f
\end{equation}
 The homomorphism $\alpha_k$ sends it to $(t_{i_1}  - 1)   ...    (t_{i_{k-1}}  - 1)    (t_{i_k}  - 1) f$.

 Let us suppose first that $a >0$, i.e. $A$ is an infinite group. 
We will  use the induction on both $m$ and $k$. The case $m=1$ is trivial:
if $(t-1)^k  f =0$ then $(t-1) \circ ... \circ (t-1) f =0$, 
even in the case   $f(t) \in  \frac{\Z[t_j]}{(t_j^{N_j} -1)}$.

Consider an element
\begin{equation} \label{rantt} 
P= \sum_j(t_{i_1(j)}  - 1)   \circ  ...  \circ (t_{i_k(j)}  - 1) f_j \in Ker \alpha_k
\end{equation} 
  
 Any element $f$ of the right hand side of (\ref{Lorant}) can be   written as 
  $
f'(t_2,...,t_m) + (t_1-1)  f''(t_1,...,t_m)$,  
 So writing $f_j = f'_j + (t_1-1)f''_j$ and setting $P': = \sum_j(t_{i_1(j)}  - 1)   \circ  ...  \circ (t_{i_k(j)}  - 1) f_j'$ we get $P =  P' +  (t_1-1) \circ Q''$.  
 
Further, let $P'_{I}$  be  the sum of all the terms of  $P'$  where non of the indices $i_l(j)$  equal to $1$.
Then $P' = P'_{I} +  (t_1-1) \circ Q'_{II}$.  The restriction of  $P'_{I}$ to the divisor $t_1=1$ in $Spec$ $\Z[A]$ coincides with the restriction of $P$. Therefore it belongs to  $Ker \alpha_k$ and thus  is zero  by the induction asumption for $(k,m-1)$. Since by the definition $P'_{I}$ does not depend on $t_1$, this implies  $P'_{I}=0$. 

    Thus $P = (t_1-1) \circ Q$ for some $Q \in  S^{k-1}I_A / S_{k-1}(A)$. Therefore $\alpha_{k-1}(Q) =0$ since $t_1-1$ is not a divisor of zero (we have assumed $a>0$). Thus  $Q=0$ by the induction assumption for $(k-1,m)$.  

Now let $A$ be a finite group. Decomposing $f(t)$ in (\ref{orant}) 
into a sum of monomials $(t_1-1)^{a_1}   ...   (t_q-1)^{a_q}$    we can write an element (\ref{orant}) as a sum
\begin{equation} \label{iiii}
P = \sum (t_1-1)^{b_1} \circ ... \circ  (t_k-1)^{b_k}   (t_{k+1}-1)^{b_{k+1}} 
  ...   (t_{k+l}-1)^{b_{k+l}} 
\end{equation}
We will treat for a moment this sum as element 
of $S^k\Z[t_1,...,t_m]$,
not $S^kI_A/S_k(A)$.  Let 
$$
\tilde \alpha_k(P):= \sum (t_1-1)^{b_1}   ...   \cdot (t_{k+l}-1)^{b_{k+l}} \in \Z[t_1,...,t_m]
$$
 be its image in $\Z[t_1,...,t_m]$
under the product map. 
If $\tilde \alpha_k(P)=0$ then clearly $P=0$ modulo the relations (\ref{rant})
with $f_i \in \Z[t_1,...,t_m]$. 

Since $\alpha_k(P) =0$, $\tilde \alpha_k(P)$ is a  
 linear combination of   monomials of type $(t_i^{N_i}-1) \cdot Q(t)$.  We are going to show that   one can find another presentation  $P'$ for the  element in $S^kI_A/S_k(A)$  given by  $P$ such that $\tilde \alpha_k(P') = \tilde \alpha_k(P) - (t_i^{N_i}-1) \cdot Q(t)$.  

Let us  factorize 
$t_i^{N_i}-1 = (t_i -1)\cdot p_i(t_i-1)$, where $p_i(u) = \sum c_j u^j$ is a   polynomial in one variable.  Notice that  $c_0 \not = 0$.  We must have a  monomial  in the sum $P$ which goes under the map $\tilde \alpha_k$ to $c_0 \cdot (t_i -1) \cdot Q(t)$. It can be written (modulo $S_k(A)$) as $c_0  \cdot (t_i -1) \circ R$.  Therefore the terms 
in (\ref{iiii}) which  are maped by   $\tilde \alpha_k$ to   $ c_{j-1}  (t_i -1)^j \cdot Q(t)$, $j>1$, can be written as     $ c_{j-1}(t_i -1)^j \circ R$. Therefore the sum of these terms in (\ref{iiii}) can be written
as $(t_i^{N_i}-1) \circ R$ and thus represent a zero element in $S^kI_A/S_k(A)$.
The proposition  is proved.

{\it Step 3}. Let us write $\beta$ as a composition  $\beta = i \circ \alpha$:
 $$
\otimes^2I^2_E \stackrel{i}{\longrightarrow}  \otimes^2I^2_E  \stackrel{\alpha}{\longrightarrow} I_E^4
$$
$$
i: A \otimes B \longmapsto A \otimes B^-; \quad \alpha: A \otimes B \longmapsto A \ast B 
$$
Let $S_E \subset \otimes^2I^2_E$ be the subgroup generated by the elements $(f) \otimes  (1-f)^-$, $f \in k(E)^*$.    
One has   
$$
 \otimes^2I^2_E = \frac{\otimes^2 k(E)^*}{k^*\otimes k(E)^*  + k(E)^* \otimes k^*},\qquad \frac{\otimes^2I^2_E}{S_E} = \frac{ K_2(k(E))}{ k^* \cdot k(E)^*}
$$
 Let $A_+ $ (resp $A_- $)   be the  coinvariants of the  involution $x \to -x$ on $E$  acting on a group $A$ functorially depending on $E$ (resp $A \otimes_{\Z}\lambda$, where $\lambda$ is  the standard $\Z$-line where the involution acts by  inverting the sign).
\begin{lemma} \label{step13}
$(Ker \alpha)_+ \subset S_E$.
\end{lemma}

{\bf Proof}.  
Let $x \in (Ker \alpha)_+$. Then $\partial (x) =0$, so $x$ defines an element of $ (H^0(E, {\cal K}_2)/K_2(k))_+$. But   this group is zero modulo $2$-torsion by lemma (\ref{3.7}).

 \begin{lemma} \label{step14}
$\Lambda^2I_E^2 \subset S_E$. 
\end{lemma}

{\bf Proof}. $A \otimes B +  B \otimes A$ belongs to the subgroup generated by the Steinberg relations $(f)  \otimes (1-f)$. Thus $A \otimes B^- +  B \otimes A^- \in S_E$. 
    The lemma above states that $B \otimes A^-  +  B^- \otimes A  \in S_E$. So  $A \otimes B^- -  B^- \otimes A  \in S_E$. 
 
Let $\alpha': S^2I^2_E \to I_E^4, \quad A \circ B \longmapsto A \ast B$.
To prove the theorem we need to show that $(Ker \alpha')_- \subset S_E$.

 Let $A_i \in I_E$. The element 
$$
< A_1,A_2,A_3,A_4>:=  
 (A_1\ast A_2) \circ (A_3\ast A_4) - (A_1\ast A_3) \circ (A_2\ast A_4)
$$
clearly belongs to  $Ker \alpha $, and  thus to $Ker \partial$. 

  \begin{proposition} \label{step15}
 $<A_1,A_2,A_3,A_4  > \in S_E$.
\end{proposition}

To prove this we need the following lemma.

\begin{lemma} \label{step16}
 If $A_i$ is a principal divisor for some $1 \leq i \leq 4$, then  \linebreak 
$<A_1,A_2,A_3,A_4> \in S_E$. In particular we get   a well defined homomorphism 
\begin{equation} \label{wdh}
<\cdot,\cdot,\cdot,\cdot>: \otimes^4J \to S^2I^2_E/S_E  
\end{equation}
\end{lemma}

{\bf Proof}. Let us show that if $B_0,B_1 \in I_E$, then $<B_0 \ast B_1  ,A_2,A_3,A_4> \in S_E$ 
(the other cases are similar).    
 We will  write $A \stackrel{S_E}{=}  B$ if $A -B \in S_E$.  It follows from  the proposition (\ref{ref})  that  for $A,B \in I_E^2$ and $X \in \Z[E(k)]$ one has 
  $(X\ast A) \circ B  \stackrel{S_E}{=} A  \circ (X\ast  B)$. 
So
$$
 (B_0\ast B_1\ast A_2) \circ (A_3\ast A_4) \stackrel{S_E}{=}    
 ( B_0\ast B_1) \circ ( A_2\ast A_3 \ast A_4) \stackrel{S_E}{=} 
 (B_0\ast B_1\ast A_3) \circ (A_2\ast A_4) 
$$
Since $(a) -(0) + (b) -(0) - ((a+b) -(0)) \in I^2_E$, we get (\ref{wdh}).
The lemma is proved.

{\bf Proof of the proposition (\ref{step15})}. Let $(A)$ be the image of $A \in I_E$ in the Jacobian. Using $(A^-) = -(A)$ and the previous lemma    we get  
$$
<A_1,A_2,A_3,A_4  >  -  <A_1^-,A_2^-,A_3^-,A_4^-  >\in S_E
$$
On the other hand according to the lemma (\ref{3.7}) 
$$
<A_1,A_2,A_3,A_4  >  +  <A_1^-,A_2^-,A_3^-,A_4^-  >  \in S_E
$$
So $2<A_1,A_2,A_3,A_4  > \in S_E$. Using lemma ({\ref{step16}) and the divisibility (by $2$) of $J(k)$ we  conclude that $<A_1,A_2,A_3,A_4  > \in S_E$.

\begin{proposition} \label{11111}
$Ker (S^2I^2_E \stackrel{\alpha'}{\longrightarrow} I^4_E)$ is generated by the elements 
$$
((x)\ast A) \circ ((-x)\ast B) -    A  \circ   B, \quad A,B \in I^2_E
$$
and $ <A_1,A_2,A_3,A_4  >$, $A_i \in I_E$.
\end{proposition}

{\bf Proof}. Let $T \subset I^2_E$ be the subgroup  generated by the elements  \linebreak $<A_1,A_2,A_3,A_4  > $. There is a surjective homomorphism   
$$
 S^4I_E   \longrightarrow  S^2(I^2_E) /T
, \quad  A_1  \otimes ... \otimes  A_4 \longmapsto ( A_1\ast A_2) \circ (A_3 \ast A_4 )
$$
 The proposition follows immediately from the proposition (\ref{schift}).

Theorem (\ref{mmmm}) follows from proposition (\ref{11111})  (\ref{schift}) and   (\ref{step16}). 

\begin{proposition} \label{qc}
 Assume that $k = \bar k$. 
Then 

a) There is an injective homomorphism of complexes
$B^*(E;3) \to B(E;3)$.

b)The quotient  $B(E;3)/B^*(E;3)$ is isomorphic   to the
Koszul complex
\begin{equation} \label{kc}
S^3J \longrightarrow S^2J  \otimes J \longrightarrow J  \otimes \Lambda^2J
\longrightarrow \Lambda^3J 
\end{equation}
In particular the complexes $B(E;3)$ and $B^*(E;3)$ are quasiisomorphic.
\end{proposition}

We will need
\begin{lemma} \label{destr}
  Suppose that 
$$
D:= m \sum_i(\{a_i\} - m \sum_{mb_i = a_i}\{b_i\}) \in I_E^4
$$
Then $D$ belongs to the subgroup generated by the elements $(f) \ast (1-f)^-$. \end{lemma}

 Let $[m]:E \to E$ be the isogeny  of multiplication by $m$.  
\begin{lemma} \label{destr1}
   Let $f,g \in k(E)^*$, $k = \bar k$. Then
$$
[m]^*\{f,g\} = m\{f,g\} \quad \mbox{in} \quad \frac{K_2(k(E)}{k^* \cdot k(E) }
$$
\end{lemma}

{\bf Proof of the lemma (\ref{destr1})}. One has an exact sequence
$$
0 \longrightarrow \frac{ H^0(E,{\cal K}_2)}{ K_2(k) + Tor(k^*,J)}\longrightarrow \frac{K_2(k(E))}{k^* \cdot k(E)  } \longrightarrow k^* \otimes J 
$$
 The operator  $[m]_*$ acts on the   left  (different from $0$) and right  groups   by multiplication by $m$ (see [BL], ch. 5,6 where this was proved rationally;  that proof works integrally). Further, $[m]_*$ has no Jordan blocks since $[m]_*[m]^* =m^2$. The lemma follows.

{\bf Proof of the lemma (\ref{destr})}. There is an isomorphism 
$$
j: \Bbb Z[E(k)]/I^4_E = \Bbb Z \oplus  J \oplus  S^2J   \oplus S^3J,
 \qquad \{a\} \longmapsto (1  , a ,   a\cdot a  , a \cdot a \cdot a)
$$
One has
\begin{equation} \label{dizr}
j(m(\{a\} - m\sum_{mb=a}\{b\})) = m((1-m^3),  (1-m^2)a,  (1-m)a \cdot a, 0)
\end{equation}
Using this we see that $\sum_i\{a_i\} \in I_E^4$. Thus $\sum_i\{a_i\} = \sum(f_j) \ast (g_j)^-$.       It is easy to see that $\beta \Bigl([m]^*\sum_i\{f_i,g_i\} - m\{f_i,g_i\}\Bigr) = D$. The lemma is proved.

{\bf Proof of the proposition (\ref{qc})}. a) We need only to show that $B^*_3(E)$ injects to $B_3(E)$. This boils down to the lemma (\ref{destr}) above, since in our definition of $R_3(E)$ we used only the distribution relations for $m=-1,2$.

b)  By definition $B_3(E)/B^*_3(E)$ is isomorphic to the quotient of 
$\Bbb Z[E(k)]/I^4_E$  modulo (the image of) the distribution relations. 

 Using the   computaion (\ref{dizr}) and divisibility of   $J(k)$ we see that the map $j$ maps the subgroup generated by the  distribution relations for any given $|m|>1$  surjectivly onto $2\Z \oplus J \oplus J^2$. Therefore
$
B_3(E)/B^*_3(E) = S^3J
$. 
Let us recall that $B_2(E)/k^* = S^2J$. So the terms of the quotient $B(E;3)/B^*(E;3)$ are the same as  in (\ref{kc}).  It is easy to see that the differentials coincide.  The proposition  is proved.

{\bf 3. Zagier's conjecture on $L(E,2)$ for modular elliptic
curves over $\Bbb Q$}. 
Let us  recall that for a curve $X$ over $\Bbb R$ one has
$H_{{\cal D}}^2(X/\Bbb R,\Bbb R(2)) =  H^1(X/\Bbb R,\Bbb R(1))$.  Let $\bar X := X \otimes \Bbb C$. The cup product with $\omega \in \Omega^1(\bar X)$ provides an isomorphism
of vector spaces over $\Bbb R$:
$$
H^1(X/\Bbb R,\Bbb R(1)) \longrightarrow H^0(\bar X, \Omega^1)^{\vee}
$$
So we will present elements of $H_{{\cal D}}^2(X/\Bbb R,\Bbb R(2))$ as
functionals on $H^0(\bar X, \Omega^1)$.

Bloch  constructed the regulator map
$$
r_{{\cal D}}: K_2(E) \longrightarrow H^2_{{\cal D}}(E,\Bbb R(2))
$$
If we represent an element of $K_2(E)$ as $\sum_i\{f_i,g_i\}$ (with all the
tame symbols vanish) then Beilinson's construction of the regulator 
looks as follows: 
$$
<r_{{\cal D}}\sum_i\{f_i,g_i\},\omega > =  \frac{1}{2 \pi i}\sum_i \int_{E(\Bbb C)}
\log|f_i|d \arg(g_i) \wedge \omega 
$$

Let 
 $f$ and $g$  be rational functions on $E$ such that
 $$
 (f)= \sum n_i(a_i), \quad  (g)= \sum m_j(b_j)
 $$
Let $\Gamma = H_1(E(\Bbb C), \Bbb Z)$.  We  may assume that  $\Gamma = \{\Z \oplus \Z\cdot \tau\} \subset \C$ and $z$ is the coordinate in $\C$.
Let us briefly recall  how  the regulator integral  $<r_{{\cal D}}\{f,g\},dz>$ is computed by means of the elliptic dilogarithm
 ([Bl1], see also [RSS]). 
 The intersection form on $\Gamma$ provides a pairing
\begin{equation} \label{myu}
(\cdot,\cdot): E(\Bbb C) \times \Gamma \longrightarrow S^1; \qquad (z,\gamma):= exp(\frac{2\pi i (z{\bar \gamma} - {\bar z}\gamma )}{\tau - {\bar \tau}})
\end{equation}
  Let 
$$
K_{2,1}(z;\tau):= \frac{(Im \tau)^2}{ \pi
}\sum_{\gamma \in \Gamma \backslash 0}\frac{(u,\gamma)}{\gamma^2{\bar \gamma}},
\qquad z = exp(2 \pi i u)
$$
 Then one has 
 $$
   \frac{1}{2 \pi i}\int_{E(\Bbb C)}\log|f| d \arg
g \wedge dz  
 =   \frac{1}{i\pi }\sum_{a,b \in E(\Bbb C)}v_a(f)v_b(g)K_{2,1}(a - b;\tau)
 $$
To prove this one may use that 
$$
\int_{E(\Bbb C)}\log|f| d i \arg g \wedge dz = - \int_{E(\Bbb C)}\log|f| d \log|g|
 \wedge dz 
$$
together with the following lemma and the fact that the Fourier
transform sends the convolution to the product.

\begin{lemma} \label{weil}
\begin {equation} \label{wei2}
\log |f(z)| = -\frac{Im \tau}{2\pi}\sum_{\gamma \in \Gamma \backslash
0}v_a(f)\frac{(z-a,\gamma)}{|\gamma|^2} +C_f
\end {equation}
where $C_f$ is a certain constant.
\end{lemma}

{\bf Proof}. One can get a
proof applying $\partial \bar \partial$ to the both parts of (\ref{wei2}). The constant
$C_f$ can be computed from the decomposition of $f$ on the product of theta functions
using the formula in s. 18 ch. VIII of [We]. It does not play any role in
our considerations since   
$
\int_{E(\Bbb C)} C_f \cdot d\log|g| \wedge \omega  =0
$
by the Stokes formula.

The relation between the Eisenstein-Kronecker series $K_{2,1}(z)$ and the
elliptic dilogarithm is the following ( [Bl1], [Z]): $
 K_{2,1}(z;\tau)= {\cal L}_{2,q}(z) - iJ_q(z)$, where  the function $J_q(z)$ is defined as follows. 
Let us  average the function $J(z):=
\log|z|\log|1-z|$  over the action of the group $\Bbb Z$ generated by the shift $z 
\longmapsto qz$  regularizing divergencies. We will get 
$$
J_q(z):= \sum_{n=0}^{\infty}J(q^nz) - \sum_{n=1}^{\infty}J(q^nz^{-1}) +
\frac{1}{3}(\log |q|)^2 \cdot  B_3(\frac{\log |z|}{\log |q|})  
$$
Here $B_3(x)$ is the third Bernoulli
polynomial. The function $J_q(z)$ is invariant under the shift $z \longmapsto qz$ and 
satisfies $J_q(z)  =  -J_q(z^{-1})$.

It follows from the main result of Beilinson in [B2], see also [SS2],  that
for a modular elliptic curve $E$ over $\Bbb Q$  there always exists an element
in $K_2(E)_{\Bbb Z}$ whose regulator gives (up to a standard nonzero factor) $L(E,2)$. 
So we get the formula 
$$
L(E,2) \sim_{\Bbb Q^*} \pi \cdot \sum_{i}\sum_{a,b \in E(\Bbb
C)}v_a(f^{(i)})v_b(g^{(i)}){\cal L}_{2,q}(a -b)  
$$

Finally, the results of the present section implies that the element \linebreak
$\sum_{i} \sum_{a,b \in E(\Bbb C)}v_a(f^{(i)})v_b(g^{(i)})\{a - b\}_3 \in I^4_{E(\bar \Bbb Q)}$  
must satisfy all the conditions of the theorem(\ref{zcc}). 

Theorem (\ref{zccc}) follows from the surjectivity of the map $\beta$ and the arguments above.

{\it  The integrality condition} ([BG], [SS]). 
Let $E$ be an elliptic curve over  $\Q$.   Choose a minimal regular model  $E_{\Z}$ of $E$
over $ \Z$. One has the exact sequence 
\begin{equation} \label {blgr}
K_2(E_{ \Z}) \longrightarrow K_2(E_{\Q}) \stackrel{\partial}{\longrightarrow} \oplus_{p  } K_1'(E_{p})
\end{equation}
The group $K_1'(E_{p}) \otimes \Bbb Q$ is not zero if and only if  $E_{p}$ has
a split multiplicative reduction with special fibre a N\'eron $N$-gon. 
In this case $K_1'(E_{p}) \otimes \Bbb Q = \Bbb Q$.

Consider an element   $\sum_i
\{f_i,g_i\} \in K_2(\Q(E))$ which has zero tame symbol at all points. It
defines an element of 
$H^0(E,{\cal K}_2)$. Suppose first   
 that the closure of the support of the divisors
$f_i,g_i$ is contained on the smooth part of $E_{\Z}$. Then Schappaher and Scholl ([SS]) proved that the image of 
this element  under the map 
$\partial$ 
is computed by the following formula:
$$
\partial (\sum_i \{f_i,g_i \} = \pm \frac{1}{3N}
\sum_{ \nu \in \Bbb Z/N\Bbb Z}d((f) \ast(g^-);\nu)B_3(\frac{\nu}{N})\cdot\Phi
$$
here $\Phi$ is a generator in $K_1'(E_{ p}) \otimes \Bbb Q$.

 In general one should extend  $\Q$ to $\Q((f_i),(g_i))$, which is the field of the definition of the divisors $(f_i), (g_i)$. After this we get precisely the condition (\ref{icond121}). 

Let us  explane why the expression
$\sum_{\nu \in \Bbb Z/(eN) \Bbb Z}d(P;\nu)B_3(\frac{\nu}{eN})$
does not depent on the field $L$.
Let $L$ and $L'$ are two extentions of $\Bbb Q_p$ such that
all the points are defined over them. We can assume that
$L\subset L'$ (by taking the composit).
Denote by $n_1=e_1 f_1$ degree of this extension. 
 Looking at   the Tate uniformization we see   
that points which intersect $\nu$-th component of
special fiber corresponding to the field $L$
intersect $e_1\nu$-th component of
special fiber corresponding to the field $L'$.

{\bf 4. Main results from the motivic point of view}. Let ${\cal M}{\cal M}_X$ be the (hypothetical) abelian category of all mixed
motivic sheaves over a regular scheme $X$ over a field $k$.  
 Let $\Bbb Q(-1)
:= h^2(P^1)$, $\Q(n):= \Q(1)^{\otimes n}$ and ${\cal H}:=
h^1(E)(1)$. 

The  motivic refinement of our results   is the following 

\begin{conjecture} \label{conjmainmot}
There exists a canonical quasiisomorphism in the derived category
$$
B(E,3)  \otimes \Bbb Q= RHom_{{\cal M}{\cal M}_k}(\Bbb Q(0), {\cal H}(1))
$$
\end{conjecture}

Let us explain  how it fits with  our  results. 
Let $\pi: E \longrightarrow Spec (k)$ be the structure morphism. 
There are the Tate sheaves $\Bbb Q(n)_E := \pi^{\ast} \Bbb
Q(n)$.
 Beilinson's description of Ext groups
between the Tate sheaves over $E$ gives us 
\begin{conjecture} \label{ext1}
$$
Ext^i_{{\cal M}{\cal M}(E)}(\Bbb Q(0)_E, \Bbb Q(2)_E) = gr^{\gamma}_2 K_{ 4-i}(E) \otimes \Bbb Q
$$
\end{conjecture}

 \begin{lemma}
\begin{equation} \label{mfor}
RHom_{{\cal M}{\cal M}_k}(\Bbb Q(0),   {\cal H}(1)) = RHom_{{\cal M}{\cal
M}_{E}}(\Bbb Q(0),\Bbb Q(2))^{ -} 
\end{equation}
\end{lemma}

Indeed, 
let $p:E \longrightarrow  Spec (k)$ be the canonical projection. Then
  we should have  the motivic Leray spectral sequence
$$
E_2^{p,q}= Ext^{p}_{{\cal M}{\cal M}_k}\Bigl(\Bbb Q(0), 
R^qp_{\ast}\Bbb Q( 2)\Bigr)  
$$
degenerating at $E_2$ and abutting to 
$Ext^{p+q}_{{\cal M}{\cal M}_{E }}\Bigl(\Bbb Q(0), \Bbb Q( 2)\Bigr)$.  
Noting that
$$
h^{0}(E )^{ -} = h^{2}(E )^{ -} = 0;    \qquad h^{ 1}
(E )^{ -} =  h^1(E) 
$$
we get (\ref{mfor}).
 Conjecture (\ref{ext1}) together with this lemma tell us that the cohomology of the  elliptic motivic  
complexes   are given by the formula
\begin{equation} \label{gipo++}
  R^iHom_{{\cal M}{\cal M}_k}(\Bbb Q(0),  {\cal H}( 1))
 \otimes \Bbb Q =
gr^{\gamma}_{ 2}K_{ 3-i}(E )^{ -}\otimes \Bbb Q
\end{equation}
 
    Conjecture (\ref{conjmainmot}) follows from theorem (\ref{za})   and    conjecture (\ref{ext1}), see [G2]. A very interesting and important is the following:

{\bf Problem}. To construct explicitely the general elliptic motivic complexes
$$
RHom_{{\cal M}{\cal M}_k}(\Q(0), Sym^n{\cal H}(m))
$$
For $m=1$ it is considered in [G2].
  
 {\bf 5. Degeneration to the nodal curve (compare with [Bl1], [DS])}.   
Let $k$ be an algebraicly closed field. Denote by ${\cal I}$ the group  
of rational functions such that $f(0) = f(\infty) =0$. Let $I_{k^*}$ be the augmentation
ideal of $\Bbb Z[k^*]$.

For an element $f \otimes g \in (1+{\cal I} )\otimes k(t)^*$  consider the Bloch map 
$$
\beta(f\otimes g):=  \sum_{x,y \in   k^*} v_x(f) v_y(g) \{y/x\} + v_{\infty}(g)\cdot ((f) + (f^-)) \in I^2_{k^*} 
$$
Set 
 \begin{equation} 
\Bbb Z[k^*] \stackrel{ \delta}{\longrightarrow}  k^* \otimes k^*, \qquad \{x\} \longmapsto (1-x) \otimes  x, \quad \{1\} \longmapsto 0
\end{equation}
Let $p: I_{k^*} \to k^*$ be the natural projection $\{x\} \to x$.
\begin{theorem}
The following  diagram is commutative:
$$
\begin{array}{ccc}
(1+{\cal I} )\otimes k(t)^*& \stackrel{ \partial}{\longrightarrow} &  I_{k^*} \otimes k^*\\
&&\\
\beta \downarrow &&\downarrow p \otimes id\\
&&\\
I^2_{k^*}&\stackrel{\delta}{\longrightarrow}&k^* \otimes k^*
\end{array}
$$
\end{theorem}

The proof is a direct calculation similar to (but   simpler then)  the  proof of theorem (\ref{commm}). 
 
Let $S(k)$ be the subgroup generated by the elements $(1-f) \otimes f$ where $f \in {\cal I} $ and the subgroup 
$(1+ {\cal I})\otimes k^* $. Set $B^*_2(k) :=  I^2_{k^*} /\beta(S(k))$.

\begin{lemma}
$Ker \beta \subset S(k)$ 
\end{lemma}

{\bf Proof}.  An easy analog of  theorem (\ref{ref})  for the nodal curve claims that the homomorphism
$$
q: k(t)^* \longrightarrow  (1+ {\cal I})\otimes k(t)^*, \quad f(t) \longmapsto  f(x) \otimes (x-1)
$$ 
is 
surjective. It is clear that $\beta \circ q$ is injective. The lemma is proved.

This lemma implies that $\bar \beta: (1+{\cal I})\otimes k(t)^* \to B^*_2(k)$  is an isomorphism. So we get a   morphism of complexes
$$
\begin{array}{ccc}
\frac{(1+{\cal I})\otimes k(t)^*}{ S(k)}& \stackrel{\partial}{\longrightarrow} &  k^*  \otimes k^*\\
&&\\
\bar \beta \downarrow &&\downarrow id\\
&&\\
B^*_2(k)&\stackrel{\delta}{\longrightarrow}&k^* \otimes k^*
\end{array}
$$
Such that $\bar \beta$ is surjective and $Ker \bar \beta  = Tor(k^*,k^*)$, similar to the theorem (\ref{zaza}). So we see that when $E$ degenerates to 
a nodal curve the complex $B^*(E/k,3)$ degenerates to  the complex $B_2^*( k ) \to k^* \otimes k^*$. 

Let $S_1(k)$ be the subgroup generated by the elements $(1-f) \otimes f$ where $f \in {\cal I} $.

 \begin{theorem}  Let $k = \bar k$. The identity map on $\Bbb Z[k^*]$ provides a homomorphism of groups
$B_2^*( k ) \to B_2(k)$. Its kernel is isomorphic to $S^2k^*$.   This map provides a 
quasiisomorphism of  complexes 
$$
\begin{array}{ccc}
B_2^*( k ) &  \longrightarrow & k^* \otimes k^*\\
\downarrow&&\downarrow\\
B_2( k ) & \longrightarrow  & k^*  \wedge k^*
\end{array}
$$ 
\end{theorem}

 One can show that  the group $B_2^*( k )$  
is isomorphic to  a  group  defined by  
S. Lichtenbaum in [Li]. 

It is known ([Lev]) that 
$$
K_2(\Bbb P^1_{\{0,\infty\}}, \{0,\infty\}) = \frac{(1+{\cal I})\otimes k(t)^*}{ S_1(k)}
$$
Using this  and the results of this subsection we can get a proof of Suslin's theorem for $k = \bar k$.

{\bf 5. On an elliptic analog of the $5$-term relation for the elliptic dilogarithm}. 
Notice that
$$
B_2(F) = Coker (\Z[ M_{0,5}(k)]   \stackrel{\partial}{ \longrightarrow} \Z[\Bbb G_m(k) ])
$$
Here $M_{0,5}$ is the configuration space of $5$ distict points on the projective line. 
So one may ask   how to present the group $B_3(E)$ in a similar form:
$$
B_3(E) \stackrel{?}{=} Coker (\Z[ X(k)]   \stackrel{\partial}{ \longrightarrow} \Z[E(k) ])
$$
where $X$ is a {\it finite dimensional} variety. In the our definition we have 
an infinite dimensional $X$ (more precisely it is an inductive limit of finite dimensional varieties). 

Here is a guess. 
Let us realize $E$ as a cubic in $P^2$. Let $p$ be a point in $P^2$ and $l_1,l_2,l_3$   any three lines through this point. Set $A_i := l_i \cap E$.
Let $\tilde R^*_3(E) \subset I^4_E$ be the subgroup generated by   the elements
\begin{equation} \label{lll} 
\{p;l_1,l_2,l_3\}:= A_1 \ast A_2^- + A_2 \ast A_3^- + A_3 \ast A_1^-
\end{equation} 
and those linear combinations of the elements $\{a\}_3 + \{-a\}_3 $ which lie in $I^4_E$. ( Probably  they belong to the subgroup generated  (\ref{lll})).
\begin{lemma}
$\tilde R^*_3(E) \subset   R^*_3(E)$.
\end{lemma}

{\bf Proof}. Let $f_i$ be  a linear homogeneous  equation of the line  $l_i$. Since the lines $l_1,l_2,l_3$ intersect in a point,  these equations are linearly dependent, so we may choose them in such a way that  
$f_1 + f_2  = f_3$. 
  Thus $f_1/f_3 \wedge f_2/f_3 = f_1/f_3 \wedge (1- (f_1/f_3))$. Applying
the map $\beta$ to this Steinberg relation and using the relations $\{a\}_3 + \{-a\}_3 =0$ we get  the element $\{p;l_1,l_2,l_3\}$.

One obviously has $\sum_{j=1}^4 (-1)^j \{p;l_1, ... \hat l_j, ..., l_4\} =0$, so we can assume that the line $l_3$ is, say,  a vertical line. 

\begin{conjecture}
 $\tilde R^*_3(E) =   R^*_3(E)$.
\end{conjecture}

\section{  $B_2(E)$,  $\theta$-functions and action of isogenies }

{\bf 1. Elliptic curves over $\Bbb C$}.  Let us represent $E$ as a quotient
$\Bbb C/{\Bbb Z}+{\Bbb Z}\tau$. 
Below $\xi$ denotes the
coordinate on $\Bbb C$. The canonical
trivialization  of  $(T^{*}_0E)^{12}$ is  $-16\prod (e_j
-e_i)(d\/\xi )^{12}$, 
which is equal to  $ \Delta (\tau ) (d\/\xi)^{12}$, where
$$
\Delta (\tau ) = (2\pi i)^{12}{\eta ^{24} (\tau )}, \quad \eta (\tau )=
q^{\frac{1}{24}}\prod_{n=1}^{\infty}(1-q^n) \quad (q=\exp (2 \pi i
\tau))
$$
 So the trivialization of  $(T^{*}_0E)^2$ is 
$ (2\pi i)^2{\eta ^4 (\tau )} (d\/\xi)^2$.

Now we will give  the analytic description of the  Deligne
pairing $[\ast ,\ast]$. 
Let $L_a$ be the line bundle corresponding to the
divisor $(a)-(0)$.  Choose
a representative $\alpha \in \Bbb C$ of $a$ ($\alpha$ is defined up to ${\Bbb
  Z}+{\Bbb Z}\tau$).  
Let us define $L_{\alpha}$ on $E$ as the quotient of the trivial line bundle
on  $\Bbb C$ under the action of ${\Bbb Z}+{\Bbb Z}\tau$~:
$1$ acts trivially and $\tau$ acts by multiplication
 by $\exp (2 \pi i \alpha )$. We identify $L_{\alpha +1 }$
with  $L_{\alpha }$ trivially and  $L_{\alpha +\tau }$ with
 $L_{\alpha}$ by multiplication  by $\exp (2 \pi i \xi )$.

The  fiber of the Poincare line bundle $[L_{\alpha},L_{\beta}]$  on
  $J\times J$  over the point $(\alpha, \beta)$ is equal to
$L_{\alpha}|_{\beta}\otimes L_{\alpha}^{-1}|_{0}$. This line bundle is  described as a quotient
 $$
\frac{\C \times \C \times \C }{(\Z \oplus \Z \tau) \oplus (\Z \oplus \Z \tau)  }
$$
$$
(\alpha, \beta, \lambda) \longmapsto (\alpha + m + n \tau, \beta + m' + n' \tau, \lambda \cdot exp(2\pi i(n\beta + n'\alpha + nn' \tau)))
$$
Let us  show that
$L_{\alpha}=L_{a}$. Consider a slight modification of the Jacobi $\theta$-function  ($z = exp(2 \pi i \xi)$):
$$
\theta (\xi) = \theta (\xi; \tau ) =  q^{1/12} z^{-\frac{1}{2}}
\prod_{j \geq  0}  (1-q^j z)\prod_{j  >  0} (1-q^j z^{-1})
$$
Then $\displaystyle
\frac{\theta (\xi -\alpha )}{\theta (\xi )\theta (\alpha )}$
is a section of $L_{\alpha}$ with  the required divisor $(a) - (0)$.

 Our recipe for the calculation of the element
$$<(a) - (0), (a) - (0)> \in  (T^{\ast}_0E)^{\otimes
2}\otimes_{k^{\ast}} [a\cdot a]$$
 leads to the expression 
$$
\frac{d\/\xi \cdot \theta '(0)}{\theta ^2(\alpha )}
( \frac{\theta ( -\alpha )}{d\/\xi \cdot \theta '(0 )\theta (\alpha
)})^{-1}=
-(\frac{d\/\xi \cdot \theta '(0)}{\theta (\alpha )})^2.$$
Notice that $\theta '(0) =2 \pi i \eta ^2(\tau )$
and the chosen analytic trivialization of  $(T^*_0E)^2$ is 
$ (2\pi i )^2{\eta (\tau )}^4 (d\/\xi)^2$. So the final answer
is 
$-\theta(\alpha)^{-2}$.

{\bf 2.  The Tate curves}.  
Let $K$ be a  field complete with respect 
to a discrete valuation $v$. Let ${\cal O}=\{a\in K|v(a)\leq 0\}$  be the ring of integers,
 $I=\{a\in K|v(a) < 0\}$ the maximal ideal,
 $k={\cal O}/I$ the residue field.

Let $q\in I$. According to Tate, the group 
$K^{*}/ q^{\Bbb Z}$ is  isomorphic to the group of
points of the elliptic curve $E_{q}$ over $K$ given by equation 
$ 
y^2+xy = x^3 +a_4x +a_6$, where
$$a_4 =-5\sum_{j\geq 1}\frac{j^3q^j}{1-q^j};
\quad a_6 =-\frac{1}{12}\sum_{j\geq 1}
\frac{(7j^5+5j^3)q^j}{1-q^j}$$
The discriminant and $j$-invariant of
this curve are:
$\Delta =q\prod_{j\geq 1}(1-q^j)^{24},
\quad j=\frac{1}{q} + 744 +196884q+\cdots$. 
The map of $K^{*}/ q^{\Bbb Z}$ to the group of
points of $E_q$ is defined by the
following expressions:
$$ 
x(u)=\sum_{j\in \Bbb Z  }\frac{q^ju}{(1-q^j u)^2}
-2\sum_{j\ge 1}\frac{jq^j}{1-q^j}; \quad 
y(u)=\sum_{j\in \Bbb Z}\frac{q^{2j}u^2}{(1-q^ju)^3}
+\sum_{j\ge 1}\frac{jq^j}{1-q^j}.
$$
The unity $1$ of $K^*$  maps to neutral element $0$ of the curve.

Define a function $T(u)$ on $K^{*}$
by the formula:
$$T(u)=\prod_{j\geq 0}(1-q^ju)
\prod_{j> 0}(1-q^ju^{-1})
$$
This function vanishes on  
$\{q^{\Bbb Z}\}$. It is quasiperiodic: $T(uq) = -u^{-1}T(u)$.

Let $a\in {\cal O}\setminus q{\cal O}$. Denote
by the same symbol its image in $E_q$. Then a
section $s$ of the bundle ${\cal O}_{E_q}
((a) -(0))$ can be represented by a function $f$
on $K^*$ such that $f(uq)= af(u)$.
Indeed, the periodic function
$f(u)T(u)(T(ua^{-1}))^{-1}$ has the required 
divisor.

 Like in the analytic case, the total space of 
the  Poincar\'e line bundle $T_{(a,b)}$ is isomorphic 
to the quotient 
$K^*\times K^*\times K^*$ modulo the action
of the group $\Z \oplus \Z$ generated by the following transformations:
$$
(a, b, \lambda )\to
(qa, b, b\lambda );\quad
(a, b, \lambda )\to
(a, qb, a\lambda ). 
$$

The corresponding group structure 
on the collection $T_{(a,\cdot)}$ is
defined by the law:
$(b_1,\lambda _1)\times 
 (b_2,\lambda _2)=(b_1b_2,
\lambda_1\lambda_2)$.

Let us calculate the expression 
$<(a)-(0),(a)-(0)>$.   
The divisor of the section
 $T(u)^{-1}(T(ua^{-1}))$ equals
$(a)-(0)$. The  regularized value 
of this expression at the divisor  $(a)-(0)$
is equal to:

$$ \frac{\displaystyle -\frac{d\,u}{u} 
(\prod_{j> 0}(1-q^j ))^2}
{\displaystyle\prod_{j\geq 0}(1-q^ja)
\prod_{j> 0}(1-q^ja^{-1})}\times
\left(\frac
{\displaystyle\prod_{j\geq 0}(1-q^ja^{-1})
\prod_{j> 0}(1-q^ja)}
{\displaystyle - \frac{d\,u}{u} 
(\prod_{j> 0}(1-q^j ))^2}  \right)^{-1}=
$$
$$
=-\left(\frac{\displaystyle \frac{d\,u}{u} 
(\prod_{j> 0}(1-q^j ))^2}
{\displaystyle a^{-\frac{1}{2}}\prod_{j\geq 0}(1-q^ja)
\prod_{j> 0}(1-q^ja^{-1})}\right)^2 . 
$$
In this calculation notice that $(1-u)\cdot u = -(u-1)$ modulo $(u-1)^2$, so the factor $1-u$ leads to $-\frac{du}{u}$.

The trivialization of $(T^*_0E)^{\otimes 12}$
is defined by the section
$$(\frac{d\,u}{u})^{12}\Delta=
\left(\left(
\frac{d\,u}{u}\right)^{2}q^{\frac{1}{6}}
(\prod_{j> 0}(1-q^j ))^4\right)^6$$
Hence the  needed expression is:
$$-\left(\frac{\displaystyle\frac{d\,u}{u} 
(\prod_{j> 0}(1-q^j ))^2}
{\displaystyle a^{-\frac{1}{2}}\prod_{j\geq 0}(1-q^ja)
\prod_{j> 0}(1-q^ja^{-1})}\right)^2
\left(\left(
\frac{d\,u}{u}\right)^{2}q^{\frac{1}{6}}
(\prod_{j> 0}(1-q^j ))^4\right)^{-1}$$
$$= - \Bigl(q^{\frac{1}{12}} 
a^{-\frac{1}{2}}T(a) \Bigr)^{-2}.
$$

{\bf 3. Calculation of the canonical  height}. a)   {\it Archimedean case}.
Let us calculate the archimedean height. 
The torsor
$T_{(a ,0)}$  is
trivialized. So the group $T_{(a,\cdot)}$
is isomorphic to the quotient of $\Bbb C\times \Bbb C^*$
with coordinates $(\beta ,\lambda )$ by 
the action of the group 
$\Bbb Z +\Bbb Z \tau$:
$$
1\colon (\beta ,\lambda ) \to
(\beta +1,\lambda ),\quad
\tau\colon (\beta ,\lambda )\to
(\beta +\tau,\lambda \times \exp (
2\pi i \alpha ))
$$

A homomorphism
$|(\cdot,\cdot )|_{\alpha}: \Bbb C\times \Bbb C^* \to 
\Bbb R^+$ which is 
invariant under the action of
$\Bbb Z +\Bbb Z \tau$ and coincides with the  norm $|\cdot|$
on $\Bbb C^*$
if $\beta =0$ is given by
$$
| (\beta, \lambda)|_{\alpha}=
\exp(-\pi i \frac{({\alpha}-\bar{\alpha})
({\beta}-\bar{\beta)}}{{\tau}-\bar{\tau}})
|\lambda |
$$
It defines a
homomorphism $T_{(\alpha ,\cdot)}  \to 
\Bbb R^+$.
 Therefore the value of our height at  
$\alpha$  equals to
$$
\log |(\alpha, 
-(\frac{ 1}{\theta (\alpha )})^2)|_{\alpha}=
2[-(\log|\theta (\alpha )| +\frac{\pi i}{2}\frac{(\alpha-
\bar{\alpha})^2}{\tau -\bar{\tau}})]
$$
It coincides with $2$ times the canonical
N\'eron height  from [Sil]. 

b)   {\it The nonarchimedean case}. 
Let $q\in I^n\setminus I^{n+1}$. The Tate curve over $Spec({\cal O})$ has a singular fiber over $Spec(k)$. The singular fiber $C = \cup C_j$ of the minimal Neron model is an $n$-gon. The set $I^j \backslash I^{j+1}$ represents the points of the curve $E$ over $Spec (K)$ whose restrictions to the singular fiber belong to the $j$-th component $C_j$ of the $n$-gon. 

Consider the map
$$\tilde{h}_v :K^*\times K^*\times K^*\to \Z;
\quad \tilde{h}_v(a, b, \lambda )=
v(q)v(\lambda )-v(a)v(b).$$
This map is invariant with respect to the action
of $\Z\oplus \Z$; therefore it detemines
a map $h_v$ from $T_{(\cdot ,\cdot) }$ to $\Z$.
The map $h_v$ is  a group homomorphism
on $T_{(a,\cdot) }$, its image is discrete,
hence its kernel is a maximal compact subgroup.

Let $v(a)=j$. The value of this map on $<(a)-(0),(a)-(0)>$ 
equals to:
$$v(q)v(\theta (a)^{-2})-v(a)v(a)=
n\cdot  (-2) \cdot(\frac{1}{12} n
-\frac{1}{2}j+v(1-a))-j^2
$$
$$
=-\frac{1}{6}n^2+jn-j^2 -2nv(1-a) = -n^2B_2(\frac{j}{n}) - 2nv(1-a)
$$
where $B_2(x) = x^2 - x + 1/6$ is the second Bernoulli polynomial. (Compare with 
the integrality condition given by the third Bernoulli polynomial).

{\bf 4. Functoriality of the groups $B_2(E)_{\Q}$ under the isogenies}. 
Let $n$ be an integer prime   to the characteristic of $  k$. Let $\lambda: E_1\to E_2$ be an isogeny of order $n$ between the elliptic curves $E_1$ and $E_2$. 
  
i) {\it Pull back }  $\lambda^*_{2 }: B_2(E_2)  \to B_2(E_1)[\frac{1}{n}] $.  Let $\hat \lambda$ be the dual isogeny $J_{E_2} \to J_{E_1}$.   The pull back of the basic extension 
$0 \to k^* \to B_2(E_1) \to S^2J_{E_1} \to 0$ (considered modulo $n$-torsion) under the homomorphism $\frac{1}{n}\hat \lambda \cdot \hat \lambda: S^2J_{E_2} \to S^2J_{E_1}$ provides  the map $\lambda^*_{2 }$:
$$
\begin{array}{ccccccccc}
0& \to &k^*&\to &B_2(E_2)&\to&S^2J_{E_2}&\to&0\\
&&&&&&&&\\
&&||&& \downarrow \lambda^*_{2 }&&\downarrow \frac{1}{n}\hat \lambda \cdot \hat \lambda&&\\
&&&&&&&&\\
0&\to&k^*&\to&B_2(E_1)[\frac{1}{n}]&\to&S^2J_{E_1}[\frac{1}{n}]&\to&0
\end{array}
$$
 
A more direct description of $\lambda^*_{2 }$ can be spelled as follows. We have  a natural morphism of $k^*$-
torsors: 
$$\lambda^{*}:[L,M]^{\otimes n}\to [\lambda^*L, \lambda^*M];
\quad
<s_1,s_2>^{\otimes n}\quad \to \quad
<\lambda^{-1}s_1,\lambda^{-1}s_2>$$
To check that it is   a map of torsors notice that
$$(<f\cdot s_1,s_2>)^{\otimes n} = f(s_2)^{  n} (<s_1,s_2>)^{\otimes n}$$
 and $f(s_2)^{  n} = (\lambda^{-1}f)(\lambda^{-1}s_2)$. 
 It is easy to see that these maps provide the pull-back map:
$$
\lambda_2^*: B_2(E_2)\to B_2(E_1)[\frac{1}{n}];\quad 
<s_1,s_2>\to
<\lambda^{-1}s_1,\lambda^{-1}s_2>^{\otimes \frac{1}{ n}}. 
$$

 \begin{theorem} \label{eqtr} 
$n\cdot \Bigl(\lambda_2^*\{\lambda(a)\}_2 - \sum_{ \gamma \in Ker \lambda}\{a + \gamma\}_2\Bigr) =0$
\end{theorem}

{\bf Proof}.   Any isogeny can be  presented as 
the composition of cyclic isogenies. 
 So we can  assume that $\lambda$ is cyclic. 
The  expression 
$$ f_n(a):= n\cdot\Bigl(\lambda_2^*\{\lambda(a)\}_2 - \sum_{ \gamma \in Ker \lambda}\{a + \gamma \}_2\Bigr) $$
 is a nonvanishing function on the noncompact curve
$E_1\setminus  Ker \lambda$.

Let $p_n: {\cal E} \to X_0(n)$ be the universal family
of elliptic curves  over the
modular curve $X_0(n)$,  
 and $\Lambda : {\cal E}\to {\cal E}$    the
universal  cyclic $n$-isogeny. 
The curves $E_1$ and $E_2$ are the fibers of  ${\cal E}$  over the points $(E_1, Ker \lambda)$ and $(E_2, E_2[n]/Ker \lambda)$ of $X_0(n)$, and $\lambda$ is the restriction
of $\Lambda$.
The construction above defines  an algebraic
function $F$ on the universal curve  ${\cal E}_1$.

Consider the punctured formal neighborhood
of  the cusp point in which the universal isogeny is
totally ramified. The $j$-invariant of
the restriction of the universal curves 
to this neighborhood has a pole at the 
cusp; hence, this restriction can be
described by the {\it Tate} curves
$E_{q}\to E_{q^n}$ ([Sil]).
 
Let us prove that $F_n$ equals $1$ for the Tate curves. We are dealing with the isogeny  $K^*/(q^n)^{\Bbb Z} \to K^*/q^{\Bbb Z}$. 
For the Tate curve we  expressed in s. 4.2  the pairing $<(a) -(0),(a) -(0)>$ in terms of the $\theta$-function 
$$
\theta_q (a) := q^{\frac{1}{12}} a^{-\frac{1}{2}}T(a) =   q^{1/12} a^{-\frac{1}{2}}
\prod_{j \geq  0}  (1-q^j a)\prod_{j  >  0} (1-q^j a^{-1}).
$$

{\bf Remark}. $\theta_q (a)$ is defined only up to a choice of sign, so only $(\theta_q (a))^2$ makes sense.


 So we need to prove the following proposition.

\begin{proposition}
\begin{equation} \label{TT}
\left(\frac{\prod_{ 0 \leq k < n} \theta_q (  t  \cdot   q^k )}{\Bigl(\prod_{0 \leq k  < n} \theta_{q^n} ( t \cdot q^{k })\Bigr)^{n} }\right)^2 =1
\end{equation}
\end{proposition}

 It shows that the restriction of the function $F_n$ to the preimage of neighborhood of the cusp point  equals $1$.  Hence this function
is equal to $1$ on all universal curve
$X_0(N)$; the function $f_n$ is the
restriction of $F_n$ to the fiber $E$;
therefore $f_n=1$.

{\bf Proof of the proposition}. The $\theta$-function  has the following property:
$$
\theta_q (t \cdot q^l) =   (-1)^k a^{-k} q^{-k^2/2}\theta_q (t)
$$
Therefore $
\prod_{0 \leq k < n } \theta_q (  t  \cdot  q^k )$ equals to 
\begin{equation} \label{raz}
  (-1)^{s_1(n)} t^{-s_1(n)} q^{- s_2(n)/2} q^{ n/12}t^{-n/2} \prod_{j \geq 0}( 1-q^{j} t   )^{n}\prod_{j  > 0}( 1-q^{j} t^{-1}   )^{n}
\end{equation}
 
Using   the definition and notations 
$$
s_1(n) := 1 + ... +(n-1) = \frac{n (n-1)}{2}; \quad s_2(n) := 1^2 + ... + (n-1)^2 = \frac{(n-1) n (2n-1)}{6}
$$
  we have     
$$
\theta_{q^n} (t \cdot  q^{k}) = q^{ n/12} t^{-1/2}  
 q^{ -k/2}
\prod_{j \geq 0}\Bigl( 1-q^{nj}\cdot  t   q^{k}\Bigr)\prod_{j  > 0}\Bigl( 1-q^{nj} \cdot t^{-1}  q^{-k}\Bigr)
$$
  On the other hand $\prod_{0 \leq k < n }\theta_{q^n} (  t \cdot  q^{k})$ is equal to 
$$
q^{\frac{n^2}{12}  -\frac{s_1(n) }{2}} t^{\frac{-n}{2}}   \cdot
\prod_{0 \leq k <n}\prod_{j\geq 0}\Bigl( 1-q^{nj+k} t  \Bigr)\prod_{0 \leq k <n}\prod_{j > 0}\Bigl( 1-q^{nj-k} t^{-1} \Bigr)= 
$$
\begin{equation} \label {dwa}
q^{\frac{n^2}{12} -  \frac{s_1(n) }{2} } t^{\frac{-n}{2}}  
 \prod_{j' \geq 0} ( 1-q^{j'} t   )\prod_{j'  > 0} ( 1-q^{j'} t^{-1} )
\end{equation}
Comparing (\ref{raz}) and (\ref{dwa}) we see that   it remains to check that 
$$
\Bigl( q^{\frac{n^2}{12} -\frac{ s_1(n)}{2} } t^{-n/2}\Bigr)^{n} = q^{- s_2(n)/2 + \frac{n}{12}} \cdot t^{-\frac{n}{2} - s_1(n)}
$$
The statement of the proposition follows.

 In particular when $\lambda [m]$ is the isogeny  of multiplication by $m$   the theorem gives  
 \begin{corollary} \label{divr}
Suppose $\bar k = k$. Then for any $a \in E(k)$ one has the "distribution relation"
\begin{equation} \label {dwadwa}
m(\{a\}_2 - \sum_{mb=a}\{b\}_2) =0
\end{equation}
\end{corollary}

{\bf Remark}. In this case one can define $[m]^*: B_2(E_2) \to B_2(E_1)[\frac{1}{m}]$ using the map $\frac{1}{m}[m] \circ \frac{1}{m}[m]$ in the diagram defining $\lambda_2$. Thus we have the factor $m$ instead of $m^2$ in the formula (\ref{dwadwa}).

ii) {\it The transfer map $\lambda_{2\ast}: B_2(E_1)[\frac{1}{n}]  \to B_2(E_2)[\frac{1}{n}] $}. Notice that the group $B_2(E)_{\Q}$ does not satisfy the descent property. Namely, if $k \subset K$ is a finite Galois extension then
$$
B_2(E/k)_{\Bbb Q} \hookrightarrow B_2(E/K)^{Gal(K/k)}_{\Bbb Q}
$$
but this inclusion is not an isomorphism because the group $S^2J(k)_{\Bbb Q}$ does not have the descent property.

 Suppose $k = \bar k$.  The transfer map $\lambda_{2\ast}$  should satisfy the projection formula
\begin{equation} \label{fedka}
\lambda^*_{2 } \circ \lambda_{2\ast } = n \cdot Id
\end{equation}
and should fit  into the following diagram, considered modulo $n$-torsion:
$$
\begin{array}{ccccccccc} \label{transfer}
0& \to &k^*&\to &B_2(E_1)&\to&S^2J_{E_1}&\to&0\\
&&&&&&&&\\
&& \downarrow m_n&& \downarrow\lambda_{2\ast }&&\downarrow \lambda \cdot \lambda&&\\
&&&&&&&&\\
0&\to&k^*&\to&B_2(E_2)&\to&S^2J_{E_2}&\to&0
\end{array}
$$
Here $m_n: x \to x^n$. This is  necessary in
 order to have (\ref{fedka}) on the subgroup $k^* \subset B_2(E)$. 

 Let us define the transfer map   as follows:
\begin{equation} \label{tr2}
\lambda_{2\ast }\{a\}_2:= \{\lambda(a)\}_2 - \Bigl(\lambda_2^*\{\lambda a\}_2 - n\{a\}_2\Bigr)
\end{equation}

{\bf Remark}.  Projection of $\lambda_2^*\{\lambda a\}_2 - n\{a\}_2$  to $S^2J$  equals
$\frac{1}{n}\hat \lambda \circ \lambda (a\cdot a) - n a\cdot  a =0$, so  $\lambda_2^*\{\lambda a\}_2 - n\{a\}_2 \in k^*$.

Let us show   that formula (\ref{tr2}) provides a transfer homomorphism.  
Suppose that $\sum\{a_i\}_2 +c=0$ in the group $B_2(E_1)$, where $c \in k^*$ (we write the group $B_2(E)$ additively).  We have to prove that  
$$
\sum_i\{\lambda(a_i)\}_2 - \Bigl(\sum_i\lambda_2^*\{\lambda (a_i)\}_2 - \sum_i n \{a_i\}_2\Bigr) +nc =0
$$ 
 By  the assumption $\sum_i\{a_i\}_2  =c^{-1} \in k^*$.  As we have shown before,
the expression in brackets always belongs to the subgroup $k^* \subset B_2(E_2)$. Therefore
$\sum_i\lambda_2^*\{\lambda (a_i)\}_2 \in k^*$.  Notice that $\lambda^*_2$ is injective modulo $n$-torsion and is the identity on the subgroup $k^* \subset B_2(E)$. Therefore modulo $n$-torsion
$\sum_i\{\lambda(a_i)\}_2 \in k^*$   and $
\sum_i\{\lambda(a_i)\}_2 = \sum_i\lambda_2^*\{\lambda (a_i)\}_2$.
 
Using proposition (\ref{eqtr}) one can easyly see that
\begin{equation} \label{tr1}
\lambda_{2\ast }\{a\}_2:= \{\lambda(a)\}_2 - \frac{1}{n}\sum_{\gamma \in Ker \lambda}\Bigl( 
n\{a+ \gamma\}_2 - n\{a\}_2\Bigr)
\end{equation}

{\bf 5. A presentation of the group $B_2(E)$ by generators and relations}. Let  ${\cal P}(a)$  (resp ${\cal P}^{'}(a)$) be the $x$-coordinate (resp.
$y$-coordinate) of the point $a \in E$ in     Tate's normal form
of an elliptic curve over an arbitrary field $k$:
$$
y^2 + a_1xy + a_3y = x^3 + a_2x^2 +a_4x + a_6
$$
\begin{proposition} \label{propo}.  a) If 
$a \not = b$ then
$$ 
<(a+b) - (0), (a+b) - (0)> \otimes <(a-b) - (0), (a-b) - (0)>$$
$$
\otimes <(a) - (0), (a) - (0)>^{-2}
\otimes <(b)-(0),(b)-(0)>^{-2} =
(\Delta^{-1/6}({\cal P}(a)-{\cal P}(b)))^{-2}
$$
b). If $a=b$ but $2a \not = 0$ then the left hand side is equal to $(\Delta^{-1/4}{\cal
P}'(a))^{-2}$. If $2a =0$ then we get $(\Delta^{-1/3}{\cal
P}''(a))^{-2}$
\end{proposition}

{\bf Proof}.  We will prove part a). Part b) is similar. 
Let $L_a$ be the line bundle corresponding to the divisor $(a)-(0)$. Evidently
$$
[L_{a +b}, L_{a +b}] \otimes [L_{a -b}, L_{a-b}]
\otimes [L_a , L_a ]^{-2}\otimes [L_b, L_b]^{-2}=k^*
$$
so
$$ 
<(a+b) - (0), (a+b) - (0)> \otimes <(a-b) - (0), (a-b) - (0)>
$$
$$
\otimes <(a) - (0), (a) - (0)>^{-2}
\otimes <(b)-(0),(b)-(0)>^{-2} \in k^* \otimes (T^{\ast}_0E)^{\otimes -4}
$$
We want to calculate this element. 
One has
$$ <(a+b) - (0), (a+b) - (0)> \otimes <(a-b) - (0), (a-b) - (0)>$$
$$\otimes <(a) - (0), (a) - (0)>^{-2}
\otimes <(b)-(0),(b)-(0)>^{-2} =
$$
$$ <(a+b) - (a), (a+b) - (a)> \otimes <(a+b)-(a),(a)-(0)>^{ 2}$$
$$\otimes <(a-b) - (a), (a-b) - (a)>\otimes <(a-b)-(a),(a)-(0)>^{ 2}$$
$$\otimes <(b)-(0),(b)-(0)>^{-2}$$

We have
\begin{equation} \label{lem}
<(a+b)-(b),(a+b)-(b)>\otimes <(a)-(0),(a)-(0)>^{-1} = 1 \in k^*
\end{equation}
Indeed, the left hand side is a regular 
function in $b$ on the elliptic curve and
so it is a constant; its value at $b=0$  is  $1$. 

Therefore the first, third and  last terms of the expression
above that we need to compute cancel thanks to (\ref{lem})
and we get:
$$<(a+b)+(a-b)-2(a),(a)-(0)>^{ 2}$$

Notice that $(a+b)+(a-b)-2(a)$ is the divisor of the function 
 $\Delta^{-1/6}({\cal P}(\xi- a)-{\cal P}(b))$.
Its value at the point $0$ is $\Delta^{-1/6}({\cal P}(a)-{\cal P}(b))$
and its generalized value at the point $a$ is the    trivialization we have chosen.
 
\begin{corollary} \label{surj}
Assume $k = \bar k$. Then the homomorphism (\ref{homo}) is surjective. 
\end{corollary}

Let us denote by $R_2(E)$ the
kernel of the homomorphism (\ref{homo}).
Then 
$$
B_2(E):= \frac{\Bbb Z[E(k) \backslash 0]}{R_2(E)}
$$

Let ${\tilde R}$ be the subgroup of $\Bbb [E(k) \backslash 0]$ generated by the  elements 
$$
\{a,b\}:= \{a+b\} + \{a-b\} - 2\{a\} -2 \{b\}   
$$
Notice that  $\{a,a\} = \{2a\} -
4\{a\}$ and $\{a,a\} - \{a,-a\} = 2(\{a\} - \{-a\})$.

\begin{lemma} \label{tors}
  For any abelian group $A$  the elements $\{a,b\}$ and $\{a\} - \{-a\}$
generate modulo $2$-torsion the kernel of the surjective homomorphism    
$\Z[A]  \longrightarrow S^2A  \quad \{a\}   \longmapsto a\cdot a$.
\end{lemma}

We will not use this fact later, so a (simple) proof is omitted. 

Consider the homomorphism ${\tilde R} \longrightarrow k^{\ast}$ defined by
the formulas
$$
\{a,b\} \longmapsto \Delta^{-1/3}({\cal P}(a) - {\cal P}(b))^2, \quad a \not = b;
\qquad \{a,a\} \longmapsto \Delta^{-1/2}({\cal P}^{'}(a))^2, \quad 2a \not = 0
$$
and $\{a,a\} \longmapsto (\Delta^{-1/3}{\cal P}^{''}(a))^2$ if $2a =0$. 
Thanks to corollary (\ref{surj}) this homomorphism 
is well defined. By definition the subgroup $R_2(E)$ is its  kernel.

{\bf Remark}. This is a homomorphism to the 
{\it multiplicative} group $k^*$ of the field $k$   defined  via the
additive structure of $k$.

{\bf 6. A remark on  the differential in the complex   $B^*(E,3)$}.    The restriction of $\delta_3$ to the subgroup $B^*_3(E)$ can be defined directly, without referring to the group $B_2(E)$ and the homomorphism $h$. A more complicated formula is  the price we pay.  
  
Set for general $a_i \in E(k)$ 
\begin{equation} \label {pfan}
\delta_3\Bigr(
(\{a_1\} - \{0\}) \ast (\{a_2\} - \{0\}) \ast
(\{a_3\} - \{0\}) \ast (\{a_4\} - \{0\})\Bigl) = 
\end{equation}
$$
\frac{[{\cal P}(a_1+a_2 ) - {\cal
P}(a_3-a_4)][{\cal P}(a_1+a_3) - {\cal P}(a_4)]
[{\cal P}(a_1+a_4) - {\cal P}(a_3)]}
{[{\cal P}(a_1 ) - {\cal
P}(a_3-a_4)][{\cal P}(a_1+a_2+a_3) - {\cal P}(a_4)]
[{\cal P}(a_1+a_2+a_4) - {\cal P}(a_3)]} \otimes  -1/2 \cdot  a_1 + ...
$$
where ... means three other terms obtained by cyclic permutation of
indices. Here  ${\cal P}(a):=  x(a)$  is  the $x$-coordinate of a point $a$.  

The expression for (\ref{pfan}) is  symmetric in $a_1,...,a_4$, which is not
obvious from the formula. 
Over $\Bbb C$ one can rewrite the right hand
side of (\ref{pfan})  in a more
symmetric way using the $\theta$-function:
$$
\frac{\theta(a_1+a_2 +a_3+a_4)
\theta(a_1+a_2)\theta(a_1+a_3)\theta(a_1+a_4)}{\theta(a_1+a_2
+a_3)\theta(a_1+a_2 +a_4)\theta(a_1+a_3+a_4)\theta(a_1)} \otimes
a_1 +...
$$
Morally the differential $\delta_3$ is given by  the ``formula''
$\{a\} \longmapsto - \frac{1}{2}  \theta(a)  \otimes a$ which,
unfortunately,   makes 
 no sence if we don't use the group $B_2(E)$.
The relation with (\ref{pfan}) is given by the classical formula
$$
{\cal P}(a) - {\cal P}(b) =
\frac{ \theta(a+b)\theta(a-b)}{\theta^2(a)\theta^2(b)}, \quad a \not = \pm b
$$

\vskip 3mm \noindent 
{\bf REFERENCES} 
\begin{itemize} 
\item[{[BT]}] Bass H., Tate J., {\it The Milnor ring of a global field} Springer Lect. Notes in Math., vol. 342 p. 389-447.
\item[{[B1]}] Beilinson A.A.: {\it Higher regulators and values of 
$L$-functions}, VINITI, 24 (1984), 181--238 (in Russian); 
English translation: J. Soviet Math. 30 (1985), 2036--2070.
\item[{[B2]}] Beilinson A.A.: {\it Higher regulators of modular curves}
Contemp. Math. 55, 1-34 (1986)
\item[{[BD2]}] Beilinson A.A., Deligne P, {\it Interpr\'etation
motivique de la conjecture de Zagier} in Symp. in Pure Math., v. 55, part
2, 1994,  
\item[{[BL]}] Beilinson A.A., Levin A.M.: {\it Elliptic polylogarithms}.
Symposium in pure mathematics, 1994, vol 55, part 2, 101-156.
\item[{[Bl1]}] Bloch S.: {\it Higher regulators, algebraic $K$- 
theory and zeta functions of elliptic curves}, Lect. Notes U.C. 
Irvine, 1977. 
\item[{[Bl2]}] Bloch S.: {\it A note on height pairings, Tamagawa
    numbers and the  Birch and Swinnerton-
Dyer conjecture} Invent. Math. 58 (1980) 65-76. 
\item[{[BG]}] Bloch S., Grayson D.: {\it $K_2$ and $L$-functions of
elliptic curves} Cont. Math. 55, 79-88 (1986)
\item[{[De]}] Deligne P. {\it La d\'eterminant de la cohomologie} Contemp.
Math. 67 93-178, (1987).
\item[{[DS]}] Dupont J., Sah S.H.: {\it Scissors congruences II}, 
J.\ Pure  Appl.\ Algebra, v. 25, (1982), 159--195.
 \item[{[G1]}] Goncharov A.B.: {\it Polylogarithms and motivic Galois
groups} Symposium in pure mathematics, 1994, vol 55, part 2, 43 -
96.
\item[{[G2]}] Goncharov A.B.: {\it Mixed elliptic motives}. To appear.
 \item[{[La]}] Lang. S.: {\it Foundations of diophantian geometry} Springer Verlag
\item[{[Lev]}] Levine M. {\it Relative Milnor K-theory}   K-theory, 6 (1992).
\item[{[Li]}] Lichtenbaum S.: {\it Groups related to  scissor congruence groups} Contemporary math., 1989, vol. 83,  151-157. 
\item[{[M]}] Milnor J.: {\it Introduction to  algebraic K-theory}  Ann. Math. Studies, Princeton, (1971).
\item[{[MT]}] Mazur B., Tate J.: {\it Biextensions and height pairings}
Arithmetic and Geometry, Birkhauser 28, v.1, (1983) 195-237.
 \item[{[RSS]}] Rappoport M., Schappacher N., Schneider P.: {\it Beilinson's conjectures on special values of $L$-functions}, Perspectives in Mathematics, vol 4, New York: Academic Press, 1988.  
 \item[{[SS]}] Schappacher N., Scholl, A.: {\it The boundary of the Eisenstein
symbol} Math. Ann. 1991, 290, p. 303-321.
 \item[{[SS2]}] Schappacher N., Scholl, A.: {\it    Beilinson's theorem on modular curves } in [RSS] 
\item[{[Sil]}] Silverman J.: {\it Advanced topics in the arithmetic      of elliptic curves} Springer Verlag, 
1994.
\item[{[S]}] Suslin A.A.: {\it $K_{3}$ of a field and Bloch's 
group}, Proceedings of the Steklov Institute of Mathematics 1991, Issue 4.
\item[{[We]}] Weil A. : {\it Elliptic functions according to Eisenstein and
Kronecker} Ergebnisse der Mathematik, 88, Springer 77.
\item[{[W]}] Wildeshaus J. : {\it On an elliptic analog of  Zagier's conjecture}
  To appear in Duke Math. Journal. 
\item[{[Za]}] Zarkhin Yu.G. {\it N\'eron pairing and quasicharacters}
Izvestia Acad. Nauk USSR 36, N2 (1972) 497-509.
\item[{[Z]}] Zagier D.: {\it The Bloch-Wigner-Ramakrishnan polylogarithm
function} Math. Ann. 286, 613-624 (1990)
\item[{[Z2]}] Zagier D.: {\it Polylogarithms, Dedekind zeta 
functions and the algebraic $K$-theory of fields}, Arithmetic Algebraic
Geometry (G.v.d.Geer, F.Oort, J.Steenbrink, eds.), Prog. Math.,
Vol 89, Birkhauser, Boston, 1991, pp. 391--430. 

\end{itemize}

A.G.: Max-Planck-Institute fur mathematik, Bonn, 53115, Germany; sasha@mpim-bonn.mpg.de

 Dept. of Mathematics, Brown University, Providence, 
RI 02912, USA.
 

A.L.: International Institute for Nonlinear Studies 
at Landau Institute
Vorobjovskoe sh. 2
Moscow 117940, Russia;  andrl@landau.ac.ru
\end{document}